\documentclass{article}

\usepackage{amsfonts}
\usepackage{amsthm}
\usepackage{amsmath}
\usepackage{mathtools}
\usepackage{hyperref}
\usepackage{xcolor}
\usepackage{tikz}
\usepackage{tikz-cd}
\usepackage{cite}
\usepackage[boxed,linesnumbered]{algorithm2e}
\usepackage{bbm}
\usepackage[capitalize]{cleveref}
\usepackage{dutchcal}
\usepackage{enumitem}

\theoremstyle{plain}
\newtheorem{theorem}{Theorem}[subsection]
\newtheorem{prop}[theorem]{Proposition}
\newtheorem{lemma}[theorem]{Lemma}
\newtheorem{cor}[theorem]{Corollary}
\newtheorem{conj}[theorem]{Conjecture}

\theoremstyle{definition}
\newtheorem{definition}[theorem]{Definition}

\theoremstyle{remark}
\newtheorem{remark}[theorem]{Remark}
\newtheorem{example}[theorem]{Example}

\crefname{cor}{Corollary}{Corollaries}

\providecommand{\keywords}[1]
{
  \small	
  \textbf{\textit{Keywords---}} #1
}

\providecommand{\msc}[1]
{
  \small	
  \textbf{\textit{MSC---}} #1
}

\newcommand{\Pp}{\mathbb{P}}

\newcommand{\Oc}{\mathcal{O}}
\newcommand{\F}{\mathbb{F}}

\newcommand{\N}{\mathbb{N}}
\newcommand{\Z}{\mathbb{Z}}

\newcommand{\Res}{\mathrm{Res}}
\newcommand{\res}{\mathrm{Res}}
\newcommand{\End}{\mathrm{End}}
\newcommand{\Endc}{\mathcal{End}}

\newcommand{\Pic}{\mathrm{Pic}}
\newcommand{\Spec}{\mathrm{Spec}}
\newcommand{\Ker}{\mathrm{Ker}}
\newcommand{\Ext}{\mathrm{Ext}}
\newcommand{\Hom}{\mathrm{Hom}}
\newcommand{\Homc}{\mathcal{Hom}}
\newcommand{\Tr}{\mathrm{Tr}}
\newcommand{\ap}{\mathfrak{a}}
\newcommand{\bp}{\mathfrak{b}}
\newcommand{\Frac}{\mathrm{Frac}}
\newcommand{\rep}{\mathrm{r\acute{e}p}}
\newcommand{\LP}{\mathrm{LP}}
\newcommand{\rank}{\mathrm{rank}}
\newcommand{\height}{\mathrm{ht}}
\newcommand{\restriction}{\mathrm{Rest}}
\newcommand{\conorm}{\mathrm{CoN}}
\newcommand{\lc}{\mathrm{lc}}
\newcommand{\piv}{\mathrm{piv}}

\title{Algebraic algorithms for vector bundles over curves}
\author{Mickaël Montessinos\thanks{Vilnius University, mickael.montessinos@mif.vu.lt}}
\date{\today}

\begin{document}
\maketitle
\begin{abstract}
    We represent vector bundles over a regular algebraic curve as pairs of lattices over the maximal orders of its function field and we give polynomial time algorithms for several tasks: computing determinants of vector bundles, kernels and images of global homomorphisms, isomorphisms between vector bundles, cohomology groups, extensions, and splitting into a direct sum of indecomposables. Most algorithms are deterministic except for computing isomorphisms when the base field is infinite. Some algorithms are only polynomial time if we may compute Hermite forms of pseudo-matrices in polynomial time. All algorithms rely exclusively on algebraic operations in function fields. For applications, we give an algorithm enumerating isomorphism classes of vector bundles on an elliptic curve, and to construct algebraic geometry codes over vector bundles. We implement all our algorithms into a SageMath package.
\end{abstract}


\msc{14Q05, 14H60, 14H05, 14-04, 94B27}
\section{Introduction} \label{sec:Intro}
Let \(k\) be a field. There is a well known equivalence between the category of regular projective curves over \(k\) and that of algebraic function fields, which we always mean to be of transcendence degree one. Furthermore, the theory of algebraic function fields itself resembles that of number fields, with the analogy deepening considerably when the field of constants is finite. Number theoretical treatments of the theory of algebraic curves are in fact not uncommon in the literature \cite{serre1975groupes,rosen2002number,stichtenoth2009algebraic}. This analogy carries over to algorithmic questions, as many algorithms first designed for number fields translate to function fields (see e.g \cite[Section 8.3]{guardia2013new}), although difficulties may arise when the characteristic is positive and when inseparable extensions are involved.

An algorithmic problem specific to function fields is the computation of the Riemann-Roch space of a divisor. Let \(K/k\) be an algebraic function field and let \(D\) be a divisor of \(K\). The Riemann-Roch space of \(D\) is the finite dimensional \(k\)-vector space
\[L(D) = \left\{a \in K^\times \mid (a) \geq -D\right\} \cup \{0\}.\]
Under the usual correspondence between divisors and line bundles (i.e locally-free sheaves of rank one), the space \(L(D)\) is isomorphic to \(H^0(X,\mathcal{L}(D))\), where \(X\) is a regular projective curve with function field \(K\). A line bundle \(\mathcal{L}\) on the curve \(X\) may then be seen as the data of a fractional ideal \(\mathcal{L}_{fi}\) of \(A_{fi}\), the integral closure of \(k[X]\) in \(K\), and a fractional ideal \(\mathcal{L}_\infty\) of \(A_\infty\), the integral closure in \(K\) of \(k(X)_\infty\), the valuation ring of the degree valuation in \(k(X)\). Then, one has \(H^0(X,\mathcal{L}) \simeq \mathcal{L}_{fi} \cap \mathcal{L}_\infty\). This observation leads to an algorithm for computing \(L(D)\) \cite{hess2002computing}.
        
The theory of divisors of function fields gives a purely algebraic counterpart to the theory of line bundles over curves, which suggests seeking an algebraic counterpart to the theory of vector bundles. While divisors themselves do not generalize well to larger rank, Weil gave a generalisation as matrix divisors \cite{weil1938generalisation} (see for instance \cite{tyurin1964classification} for a more modern account). Our representation for vector bundles is similar to matrix divisors, but tailored for our purposes of algorithmic manipulation.

Line bundles may be described by giving generators of their stalks at the closed points of \(X\). Such a family of generators is represented by an idèle of \(K\). For vector bundles of higher rank, this generalizes to invertible matrices with coefficients in the ring of adèles of \(K\). This is the approach followed by \cite{weng2018codes,weng2018adelic}. One may also consider pairs of lattices over \(A_{fi}\) and \(A_\infty\), where by lattices we mean full rank projective sub-modules of \(K^r\). This analogy leads to a natural generalisation of Hess' algorithm for computing the space of global sections of a vector bundle, and was in fact used in \cite[Lemma 25]{ivanyos2018computing}, where the vector bundle is a maximal order in a \(K\)-algebra isomorphic to \(M_n(K)\). 

In this work, we present a computational representation of vector bundles as pairs of lattices. We also use the adelic framework discussed above for proving theoretical results. After introducing preliminary notions and quoting useful results in \Cref{sec:Preliminaries}, theoretical results in an adelic framework are stated and proved in \Cref{sec:VectorBundles}. A practical algorithmic representation, along with several algorithms for it, are discussed in \Cref{sec:Algorithms}. Finally, in \Cref{sec:Applications} we give explicit examples of some constructions mentioned in \Cref{sec:IntroApp}

\subsection{Our contributions}
    We give a computational representation of vector bundles over regular projective curves, as well as polynomial time algorithms for the following tasks:
    \begin{itemize}
        \item Compute the determinant and degree of a vector bundle. (\Cref{thm:BasicAlgos,rem:AlgosEasy})
        \item Compute dual bundles, bundles of homomorphisms, direct sums and tensor products. (\Cref{thm:BasicAlgos,rem:AlgosEasy,rem:AlgoHom})
        \item Compute direct and inverse images of vector bundles by a separable morphism of algebraic curves. (\Cref{sec:RestrictionConormVB,sec:DirectImages}, in particular \Cref{thm:ResConorm})
            \item Compute the cohomology groups of a vector bundle and compute the Serre duality isomorphism. We stress that our algorithm for computing Serre Duality is only valid when the base field \(k\) is finite. (\Cref{cor:H0,thm:H1})
            \item Compute extensions of vector bundles from their cohomological representation. Since this relies on the computation of \(H^1\), it is also only valid when \(k\) is finite. (\Cref{cor:ExtensionsAlgorithm})
            \item If the constant base field \(k\) is infinite (or large enough), test if two vector bundles are isomorphic and compute an isomorphism if it exists. This algorithm is probabilistic of Monte Carlo type. However, if it outputs an isomorphism, it is always correct. (\Cref{thm:AlgoIsomLargeField})
    \end{itemize}
            Some algorithmic tasks require computing Hermite normal forms for matrices over maximal orders in algebraic function fields. We argue that it is likely that this task may be done in polynomial time but we leave this problem for further work. This is discussed in \Cref{sec:Hermite} and stated as \Cref{conj:Hermite}. The following tasks may be done in polynomial time, given access to an oracle for computing Hermite forms of matrices:
    \begin{itemize}
        \item Compute the image and kernel of a homomorphism of vector bundles. (\Cref{cor:AlgoImage,cor:AlgoKernel})
        \item If the base field \(k\) is finite and \(L\) is a vector bundle, compute an isomorphism \(L \simeq \bigoplus L_i^{n_i}\), where the \(L_i\) are indecomposable vector bundles. (\Cref{thm:Splitting})
        \item If the base field \(k\) is finite, test that two vector bundles are isomorphic and compute an isomorphism if it exists. This algorithm is deterministic. (\Cref{cor:IsomFiniteField})
    \end{itemize}

We also provide an implementation of our algorithms as a package for SageMath \cite{sagemath}\footnote{\url{https://git.disroot.org/montessiel/vector-bundles-sagemath}}.

\subsection{Related work}
In a more general scope, computations on coherent sheaves over projective varieties may be achieved using their description as graded modules \cite{serre1955faisceaux}. Increasingly efficient methods have been developped to compute the cohomology of sheaves, see for instance \cite{smith2000computing, decker2002sheaf, motsak2010graded}.

Our approach being smaller in scope, we can afford a more straightforward and explicit representation of vector bundles, with direct ties into the description of their isomorphism classes as the double quotient
   \[GL_n(K)\backslash GL_n(\mathbb{A}) / GL_n(\mathcal{O}_\mathbb{A}).\]

An algorithm for deciding if a coherent sheaf is semi-stable is introduced in \cite{brenner2022deciding}. While tailored for the usual computation paradigm discussed above, the algorithm relies mostly on computing pullbacks of bundles and computing twists and global sections, which are tasks accessible in our setting. Over finite fields, the algorithm may conclude that a sheaf is semi-stable or that it is not strongly semi-stable (as defined by the authors). More work is required to give a full decision algorithm for schemes over finite fields, and to allow for the computation of a destabilizing subsheaf.

   \subsection{Applications}\label{sec:IntroApp}
A long line of work has been dedicated to the geometric study of moduli spaces and moduli stacks of vector bundles over curves and more generally over algebraic varieties. While efforts are usually focused on the case of varieties over the field of complex numbers, curves over finite fields were famously used in \cite{harder1975cohomology} to compute Betti numbers of moduli spaces over complex curves. An algorithmic treatment of vector bundles may be used to give algorithmic enumerations of isomorphism classes of vector bundles on a curve over a finite field. In fact, the material present in this article is sufficient to handle the case of genus \(0\), where indecomposable vector bundles are simply the line bundles by a famous theorem of Grothendieck, and the case of genus \(1\), for which Atiyah gave an explicit description \cite{atiyah1957vector}. The case of genus \(1\) is discussed in detail in \Cref{sec:EllipticCurves}.

Higher rank vector bundles on algebraic curves over finite fields can be used to construct error correction codes, generalizing the usual algebraic geometry codes. The idea was first suggested in \cite{savin2008algebraic} and then further developed in works such as \cite{ballico2008vector, nakashima2010ag}, with an independent discovery in \cite{weng2018codes}. Our algorithms allow for the construction of so-called weakly-stable vector bundles that are useful for these purposes, and to compute spaces of global sections necessary for encoding and decoding words. We give an example of such a construction in \Cref{sec:Codes}.

Finally, we mention the application of vector bundles for studying the explicit isomorphism problem over function fields. The problem is, given a \(K\)-algebra \(A\) isomorphic to a matrix algebra, to find an explicit isomorphism \(\varphi \colon A \to M_r(K)\). In \cite{ivanyos2018computing}, the authors give a polynomial time algorithm for the case \(K = \F_q(x)\). Their approach may be rephrased as follows:
\begin{enumerate}
    \item Compute a maximal order over \(\Pp_{\F_q}^1\) of \(K\)-algebra \(A\) of dimension \(r^2\).
    \item This is the sheaf of endomorphisms of some vector bundle of rank \(r\) over \(\Pp_{\F_q}^1\).
    \item Such a vector bundle splits as a direct sum of line bundles by Grothendieck's theorem.
    \item Compute the \(\F_q\)-algebra of global endomorphisms of the order. It must contain a rank one idempotent of \(A\).
    \item Use this idempotent element to compute an explicit isomorphism \(A \to M_r(K)\).
\end{enumerate}
In the case of a field of higher genus, Grothendieck's theorem is no longer true and the algebra of global sections of a maximal order may be isomorphic to \(\F_q\), for instance in the case of the sheaf of endomorphisms of a stable vector bundle. , In this case, the 4th step of the strategy above fails. A possible approach to overcome this obstacle is to find a way to, given a maximal order that is the sheaf of endomorphisms of an indecomposable vector bundle, compute a different maximal order which will have global zero divisors. It is the hope of the author that a better understanding of the algorithmic theory of vector bundles will yield helpful results in this direction.

\subsection{Further work}
    For algebraic curves of genus larger than \(1\), the early literature on moduli spaces of vector bundles focuses on constructing semi-stable and stable bundles. If one wishes to include all vector bundles over the curve, a moduli stack must be used instead. In further work, we wish to investigate the algorithmic construction and enumeration of stable vector bundles over a curve of genus larger than \(1\), and the algorithmic detection of (semi-)stability for a vector bundle.

    As discussed above, we hope to be able to use this setting for investigating the algorithmic behavior of maximal orders in matrix algebras, with applications to the explicit isomorphism problem.

    Finally, we wish to give a polynomial time algorithm for computing Hermite normal forms of matrices over maximal orders of function fields and prove \Cref{conj:Hermite}.

\subsection{Setting and Notations}
    Let \(k\) be a perfect field. We let \(k(x)\) be the rational function field over \(k\) and \(K \coloneq k(x,y)\) be a separable extension of degree \(n\) of \(k(x)\). We let \(k[x]\) be the ring of polynomials over \(k\) and \(k(x)_\infty\) be the valuation ring of the \(-\deg\) valuation in \(k(x)\). Then, we denote by \(A_{fi}\) the integral closure of \(k[x]\) in \(K\) and by \(A_\infty\) the integral closure of \(k(x)_\infty\).

    The set of places (see e.g \cite[Definition 1.1.8]{stichtenoth2009algebraic}) of any function field \(L\) is denoted by \(M_L\), or simply \(M\) when the field is clear from context. Then, \(M_L^{fi}\) is the set of finite places and \(M_L^\infty\) is the set of infinite places. For any place \(P \in M_L\), \(k_P\) is the residue field at \(P\), \(v_P\) is the normalized valuation at \(P\) (meaning \(v_P(K^\times) = \Z\)), and \(A_P\) is the valuation ring of \(P\). The degree of \(P\) is \(\deg P = [k_P:k]\). In general, we assume a choice of local uniformizer \(\pi_P\) for each place \(P \in K\). This is consistent with the computational situation where one has a deterministic algorithm for computing a local uniformizer at any place.

    While most of our treatment is focused on algebra in the field \(K\), the objects we handle represent vector bundles over a normal projective curve \(X\) over \(k\) with function field \(K\). The presentation of \(K\) as the separable extension \(k(x,y)\) then yields a separable map \(x\colon X \to \Pp_k^1\).

    For the remainder of this article, whenever we mention fields \(k\) and \(K\), we mean fields defined as above unless we specify otherwise.

    If \(g = (g_{ij})\) and \(g' = (g'_{ij})\) are matrices and \(g\) has size \(r_1 \times r_2\), we write \(g \otimes g'\) for the Kronecker product
    \[g \otimes g' = \begin{pmatrix} g_{11} g' & \hdots & g_{1r_2} g' \\
               \vdots & \ddots & \vdots \\
               g_{r1} g' & \hdots & g_{r_1 r_2} g'
           \end{pmatrix}\]
    and we set \(g \oplus g'\) to be the diagonal join
    \[
        g \oplus g' = \begin{pmatrix} g & (0) \\ (0) & g'\end{pmatrix}.\]

\paragraph{Acknowledgments}
The author wishes to thank Lin Weng for his encouragments and for conversations about his work on adelic matrices and vector bundles, and Holger Brenner for sharing his work on the detection of semi-stability.

\section{Preliminaries} \label{sec:Preliminaries}
    \subsection{Algorithms for function fields}\label{sec:BasicAlgos}
    All our algorithms rely on algebraic computations in an algebraic function field over a perfect base field \(k\). We further assume that the elements of \(k\) may be represented exactly, and that we have efficient algorithms for computing arithmetic operations in \(k\) and for factoring polynomials in \(k[x]\). The complexity of our algorithms is expressed in terms of number of arithmetic operations in \(k\) and polynomial factorizations in \(k[x]\). In particular, our algorithms have polynomial complexity in the usual sense if \(k\) is a finite field.

    The field \(K\) is represented computationally as the vector space \(k(x)^n\), using the basis \(1,y,\hdots,y^{n-1}\). A fractional ideal of \(A_{fi}\) or \(A_\infty\) is represented by the Hermite normal form of its basis over \(k[x]\) or \(k(x)_\infty\).

    While several strategies may be used for most computations in a function field, we simply record that the following problems may be solved in polynomial time:
    \begin{itemize}
        \item Compute a fractional ideal of \(A_{fi}\) or \(A_\infty\) given a finite set of generators.
        \item Compute the divisor associated to a fractional ideal of \(A_{fi}\). 
        \item Efficient Chinese Remainder Theorem: given places \(P_1,\hdots,P_r \in M\), integers \(n_1,\hdots,n_r\) and elements \(a_1,\hdots,a_r \in K^\times\), find \(a \in K^\times\) such that \(v_{P_i}(a-a_1) \geq n_i\) for all \(1 \leq i \leq r\).
        \item Given a place \(P \in M\), a local uniformizer \(\pi_P\), an integer \(i \in \Z\) and some \(a \in K^\times\), compute the coefficient \(a^{(i)}\) of the formal power series \(a = \sum_{i \geq v_P(a)} a^{(i)} \pi_P^i\).
    \end{itemize}

    The first three problems admit natural counterparts in the algorithmic theory of number fields. Algorithms from \cite{cohen1993course, cohen2000advanced} can then be adapted to the function field case. We also note the more recent work \cite{guardia2013new} which is applicable when \(k\) is a finite field. In both cases, the algorithms presented for number fields become fully polynomial over function fields since we assume polynomial time algorithms for factoring polynomials. Computing coefficients of power series may be done for instance using \cite[Lemma 9 and Algorithm 27]{hess2002algorithm}.

    \subsection{Hermite normal form}\label{sec:Hermite}
    The algorithms given in \cite{cohen2000advanced} and mentioned in \Cref{sec:BasicAlgos} require computing the Hermite normal form of matrices with coefficients in \(k[x]\) or in the valuation ring \(k(x)_\infty\). The Hermite normal form of a matrix with coefficients in \(k[x]\) may be computed in polynomial time using \cite[Theorem 5]{gupta2012triangular}. The ring \(k(x)_\infty\) is a \(DVR\), so the Hermite normal form may also easily be computed over this ring.

    The algorithms given in \Cref{sec:AlgosHomomorphisms} will rely of the computation of Hermite normal forms of matrices and pseudo-matrices over the Dedekind domains \(A_{fi}\) and \(A_\infty\). We briefly recall the relevant definitions and results, and refer the reader to \cite[Sections 1.4 and 1.5]{cohen2000advanced} for details. Note that we include the case of matrices not of full rank.

    For the rest of this section, \(A\) is a Dedekind domain with fraction field \(K\). We first give a definition of pseudo-matrices:
    \begin{definition}[{\cite[Definition 1.4.5]{cohen2000advanced}}]
        \begin{enumerate}
            \item A pseudo-matrix of size \(r \times n\) over \(A\) is a pair \((\ap,M)\), where \(\ap = (\ap_j)_{j=1}^n\) are fractional ideals of \(A\) and \(M \in M_{n,r}(K)\).
            \item The map associated with such a pseudo-matrix is the map \(f\colon \ap_1,\hdots,\ap_n \to K^r\) defined by \(f(a_1,\hdots,a_n) = \sum_{j=1}^n a_j M_j\), where the \(M_j\) are the columns of \(M\).
            \item The module associated with this pseudo-matrix is the module
                \[L = \bigoplus_{j=1}^r \ap_j M_j\]
                it is the image of the map \(f\) in \(K^n\).
            \item The kernel of the pseudo-matrix \((\ap,M)\) is the kernel of the map \(f\).
        \end{enumerate}
    \end{definition}

    \begin{remark}
        If \(A\) is a PID, we may always turn a pseudo-matrix into a matrix by computing generators of its coefficient ideals. We state the results in a unified manner but we always assume that the coefficient ideals are trivial when the base ring is \(A_\infty\).
    \end{remark}

    A detailed definition of the Hermite normal form will not be needed, as we will only use the following known facts:
    \begin{prop}\label{prop:Hermite}
        \begin{enumerate}
            \item Let \((\ap,M),(\ap',M')\) be two pseudo-matrices over \(A\). Then the modules generated by these pseudo-matrices are equal if and only if their Hermite normal forms are equal as well \cite[1.5.2 (2)]{cohen2000advanced}.
            \item The image and kernel of a pseudo-matrix may be computed in polynomial time from its Hermite normal form. \cite[1.5.2 (5)]{cohen2000advanced}
        \end{enumerate}
    \end{prop}

    We now discuss the computation of Hermite normal forms for pseudo-matrices over the rings \(A_{fi}\) and \(A_\infty\). The analogous problem for maximal orders of number fields was conjectured in \cite{cohen1996hermite} to be feasible in polynomial time, and it was proved in \cite{biasse2017computation}. In the case of maximal orders in a function field, while no such result exists in the literature to the best of our knowledge, it seems likely that the methods used for the number field case may be adapted, using the computation of Popov forms of matrices over \(k[x]\) (see \Cref{sec:Popov}) to replace the use of LLL reduction. We leave this investigation to further work. For now, we state:

    \begin{conj}\label{conj:Hermite}
        There exists a polynomial-time algorithm to compute the Hermite normal form of a pseudo-matrix over \(A_{fi}\).
    \end{conj}

    \subsection{Popov form of matrices of polynomials}\label{sec:Popov}
    In order to relate lattices over the rings \(A_{fi}\) and \(A_\infty\), we need to compute bases that are orthogonal in some sense. While LLL reduction would be used over number fields, the function field equivalent we use here is the Popov form of matrices. We follow the exposition given in \cite{sarkar2011normalization}, except all statements are transposed since our convention will be to have the columns of matrices represent basis elements.

    \begin{definition} \label{def:NormVector} Let \(v = (v_i) \in M_{r,1}(k[x])\) be a column vector. We define the following:
        \begin{itemize}
            \item The \emph{norm} of \(v\) as \(|v| = \max_{i=1}^r \deg(v_i)\).
            \item The vector \(\lc(v) \in M_{r,1}(k)\) is the vector whose \(i\)-th entry is the coefficient of degree \(|v|\) of the polynomial \(v_i\).
            \item The \emph{pivot index} of vector \(v\), denoted by \(\piv(v)\) is the largest \(i\) such that \(\deg v_i = |v|\).
        \end{itemize}
    \end{definition}

    \begin{definition}[{\cite[Definition 2]{sarkar2011normalization}}]
        Let \(M \in M_r(k[x])\). We say that the matrix \(M\) is reduced if the matrix
        \[\lc(M) = \begin{pmatrix} \lc(v_1) & \hdots & \lc(v_r) \end{pmatrix}\]
        is invertible.
        Let \(v_1,\hdots,v_r\) be the columns of \(M\). We say that the matrix \(M\) is in Popov form if it is reduced and the following conditions are satisfied:
        \begin{enumerate}
            \item The pivot indices \(\piv(v_1),\hdots,\piv(v_r)\) are distinct.
            \item The pivot entries \(v_{i,\piv(v_i)}\) are monic.
            \item For \(1 \leq i < r\), \(|v_i| \leq |v_{i+1}|\), and if \(|v_i| = |v_{i+1}|\), then \(\piv(v_i) < \piv(v_{i+1})\).
            \item The entries of \(v\) that are not the pivot of their column have degree lesser than the entry of the same row which is the pivot of its own column.
        \end{enumerate}
    \end{definition}

    The reason reduced matrices are relevant to us is the following statement, which is a form of orthogonality:
    \begin{prop}[{\cite[Theorem 6.3-13]{kailath1980linear}}]\label{prop:DegreePredictability}
        Let \(v_1,\hdots,v_r\) be the columns of a reduced matrix. Let \(a_1,\hdots,a_r \in k[x]\). Then
        \[\left|\sum_{i=1}^r a_i v_i \right| = \max_{1 \leq i \leq r}\left(\deg(a_1) + |v_i|\right).\]
    \end{prop}

    Computing reduced matrices would be sufficient for most applications, but a matrix may be right-equivalent to several different reduced matrices. Computing the Popov form of a matrix instead ensures unicity, and may be desirable in a computational context.
    \begin{prop}
        Let \(M \in M_r(k[x])\) be nonsingular. Then there exist unique matrices \(U \in GL_r(k[x])\) and \(P \in M_r(k[x])\) such that \(P = MU\) and the matrix \(P\) is in Popov form.
    \end{prop}

    Reduced and Popov forms of matrices may be computed efficiently. In the following, the notation \(\tilde{O}\) means we omit logarithmic factors, \(\omega\) is the exponent of the cost of matrix multiplicaiton in \(k\), \(M(d)\) is the cost of multiplication of two polynomials of degree at most \(d\) and \(B(d)\) is the cost of an extended gcd computation for two polynomials of degree at most \(d\).
    \begin{prop}
        Let \(M \in M_r(k[x])\) with entries of degree no larger than \(d \in \N\).
        \begin{itemize}
            \item A reduced matrix right-equivalent to \(M\) may be computed in \(\tilde{O}(n^\omega (M(d) + B(d)))\) operations in \(k\). \cite{gupta2012triangular}
            \item If \(M\) is reduced, the Popov form of \(M\) maybe computed in \(\tilde{O}(n^\omega d)\) operations in \(k\). \cite{sarkar2011normalization}
        \end{itemize}
    \end{prop}
    \subsection{Heights}
    Let \(a \in K^\times\). We define the \emph{height} of \(a\) as \[\height_K(a) = \sum_{P \in M} \max(v_P(a),0) = \sum_{P \in M} -\min(v_P(a),0).\] We simply write \(\height(a)\) if the field \(K\) is clear from context. Observe that for a finite separable extension \(K'/K\), if \(a \in K^\times\), \(\height_{K'}(a) = [K':K]\height_K(a)\).

    While analogous to the notion of heights in number fields, this is more usually called \emph{degree} in the literature on function fields. In this work, we already consider degrees of vector bundles in various forms and connect it to degrees of divisors and répartitions. Therefore, we use the term height for this notion.
    
    We observe readily that for \(a,b \in K^\times\), \(\height(a+b) \leq \height(a) + \height(b)\) and \(\height(ab) \leq \height(a) + \height(b)\).

    In the case \(K = k(x)\), let \(r = \frac{p}{q} \in k(x)\) with \(p\) and \(q\) coprime non-zero polynomials in \(k[x]\). Then \(\height(r) = \deg(p) + \deg(q)\) and therefore the height of \(a\) is linear in the size of the representation of \(r\).

    Then, we connect the height of an element \(a \in K^\times\) with the size of its computational representation:
    \begin{prop}\label{prop:BoundedHeight}
        Let \(a \in K^\times\). Then \(\height(a)\) is polynomial in the size of the representation of the field \(K\) and of the function \(a\).
    \end{prop}

    \begin{proof}
        First, we compute \(\height(y)\). Let \(\chi_y = \sum_{i=0}^n c_i T^i\) be the minimal polynomial of \(y\) over \(k(x)\). Observe that \([K:k(y)] \leq \max_i \height(c_i)\). Thus, by \cite[Theorem 1.4.11]{stichtenoth2009algebraic}, \(\height_K(y) \leq \max_i \height(c_i)\). The size of the representation of the field \(K\) is bounded by \(\sum_{i=0}^{n} \height(c_i)\), so \(\height(y)\) is bounded by the size of the representation of \(K\).
        Let \(a = \sum_{i=0}^{n-1} a_i y^i\). Then we compute
        \begin{align*}
            \height(a) &\leq \sum_{i=0}^{n-1} \height_K(a_i) + i \height(y) \\
                       &\leq \sum_{i=0}^{n-1}(\height(a_i)) + \frac{n(n-1)}{2} \height(y).
        \end{align*}
    \end{proof}

    We will need the following lemma in the sequel:
    \begin{lemma}\label{lemma:HeightOfDegree}
        Let \(g \in GL_r(K)\). Then \(\height(\det(g))\) is polynomial in the size of representation of \(g\).
    \end{lemma}

    \begin{proof}
        A representation of the determinant of \(g\) in the basis \(1,y,\hdots,y^{n-1}\) may be computed in polynomial time, and therefore has polynomial size. The result then follows from \cref{prop:BoundedHeight}.
    \end{proof}

    \subsection{Répartitions}\label{sec:Repartitions}
    The words répartition and adèles are sometimes used interchangeably in the literature. When working over a global field, they always mean to take a restricted product over the set of places of the field. Most usually, this is the product of the completions of the field, but one may also work with mere copies of the field. Since our work is computational in nature, we avoid taking completions to preseve exact representations. In this work, we use the term \emph{répartition} to emphasize this. References for répartitions in algebraic function fields are \cite{serre1975groupes,stichtenoth2009algebraic}.

    \begin{definition}
	    The ring \(R_K\) of répartitions of \(K\) (simply written \(R\) when there is no ambiguity) is the restricted product
	    \[
		    R = \tilde{\prod_{P \in M}} K
	    ,\]
	    where the restriction means that for an element \((r_P) \in R\), all but finitely many of the \(r_P\) lie in \(A_P\).

	   The subring of \(R\) of integral répartitions is the product
	   \[
		   A_R = \prod_{P \in M} A_P.
	   \]
    \end{definition}

    Now, we may define the degree of a répartition
    \begin{definition}
	    Let \(r \in R^\times\), we set
	    \[\deg(r) = \sum_{P \in M} v_P(r_P) \deg(P).\]
    \end{definition}
    The degree is well defined since, for all but finitely many \(P \in M\), an invertible répartition \(r \in R^\times\) lies in \(A_P^\times\) and therefore has valuation zero at \(P\).

    In \cref{sec:VectorBundles}, we will give an explicit statement of Serre duality for vector bundle in terms of répartitions. For this purpose, we will need to recall some definitions relating répartitions and differentials.
    \begin{definition}
	    Let \(\omega\) be a differential of \(K\). For any \(P \in M\), let \(f_P \in K\) be defined by the relation \(\omega = f_P d\pi_P\). Then, \((f_P)_{P \in X}\) is the \emph{répartition of \(\omega\)}, which we denote by \(\iota(\omega)\).
    \end{definition}

    Of course, the répartition \(\iota(\omega)\) depends on the choice of local uniformizers. However, such a choice leaves the residue of a differential invariant. Hence, we may define:
    \begin{definition}
	    Let \(r \in R\). The \emph{residue} of \(r\) is defined as the sum
	    \[\res(r) = \sum_{P \in M} \Tr_{k_P/k}(\res_{P}(r_P)),\]
        where \(\res_P(r_P)\) is the trace over \(k\) of the coefficient of degree \(-1\) in the formal series in \(\pi_P\) representing \(r_P\).
    \end{definition}

    \begin{remark}
	    Our definition of the residue of a r\'epartition mimics that of the residue of a differential. The residue of a répartition is always well defined as a répartition must have a nonnegative valuation at all but finitely many places. It is immediate by the residue theorem that if \(\omega\) is a differential, \(\res(\iota(\omega)) = 0\).
    \end{remark}

    In \Cref{sec:H1}, we will need to compute répartitions with prescribed residues. Our strategy will be to focus on an infinite place of \(K\). We introduce the following useful notation:
    \begin{definition}
        An \emph{infinite répartition} is a répartition \(r \in R\) such that \(r_P = 0\) for all \(P \in M^{fi}\) and there is an \(a \in K\) such that \(r_P = a\) for all \(P \in M^\infty\). In this case, we denote the r\'epartition \(r\) by \(a_\infty\). We also extend this notation the obvious way to vectors with coefficients in \(R\).
    \end{definition}

    In order to use répartitions to describe vector bundles of rank larger than \(1\), we use matrices taking coefficients in \(R\). Such techniques were already discussed in \cite{sugahara2012adelic,weng2018zeta,weng2018codes,weng2018adelic} for matrices with coefficients in the ring of adèles. The properties we need to establish are often analogous to some results from the references above, but we give our own proofs for completeness and to account for the change from adèles to répartitions.

    A matrix \(M \in M_{r_1,r_2}(R)\) may be seen as a family \((M_P)_{P \in M}\) of matrices in \(M_{r_1,r_2}(K)\) with the extra condition that at most finitely many of the \(M_P\) do not lie in \(M_{r_1,r_2}(A_P)\).

    We take note of an easy lemma
    \begin{lemma}\label{lemma:GLnR}
	    A matrix \(M \in M_n(R)\) is invertible if and only if it lies in \(\prod_{P \in M} GL_n(K)\) and all but finitely many of the \(M_P\) lie in \(GL_n(A_P)\).
    \end{lemma}

    \begin{proof}
	    The determinant \(d = \det M\) is invertible in \(R\) if and only in it lies in \(\prod_{P \in M} K^\times\) and for all but finitely many \(P\), \(d_P \in A_P^\times\). The result follows readily.
    \end{proof}

\section{Vector bundles and lattices}\label{sec:VectorBundles}
We begin with a definition of four equivalent categories:

\begin{definition}\label{def:VBCategories}
	\begin{enumerate}
		\item A \emph{vector bundle} over \(X\) is a locally free sheaf of coherent \(\mathcal{O}_X\)-modules. A map of vector bundles is simply a homomorphism of \(\Oc_X\)-modules. \label{def:VectorBundle}
		\item Let \(K_X\) be the constant sheaf equal to \(K\) over \(X\). Then, an \emph{\(\Oc_X\)-lattice} is a subsheaf of \(K_X^r\) that is a locally free \(\Oc_X\)-module of rank \(r\), for some \(r \in \N\). A map between \(\Oc_X\)-lattices \(\mathcal{L_1}\) and \(\mathcal{L_2}\) of ranks \(r_1\) and \(r_2\) is a \(\Oc_X\)-module homomorphism \(f \colon K_X^{r_1} \to K_X^{r_2}\) such that \(f(\mathcal{L}_1) \subset \mathcal{L}_2\). \label{def:OxLattice}
        \item An \emph{\(A_R\)-lattice} is a free \(A_R\)-submodule \((L_P)_{P \in M}\) of \(R^r\) for some \(r\). A map between \(A_R\)-lattices \(L_1\) and \(L_2\) of ranks \(r_1\) and \(r_2\) is a map \(f\colon K^{r_1} \to K^{r_2}\) such that \(f(L_1) \subset L_2\) when \(f\) is extended to \(R^{r_1}\) by pointwise application. \label{def:ARLattice}
        \item A \emph{lattice pair} \(L\) of \(K\) is the data of an \(A_{fi}\)-lattice \(L_{fi}\) and an \(A_\infty\)-lattice \(L_\infty\) of equal ranks. A map between lattice pairs \(L\) of rank \(r\) and \(L'\) of rank \(r'\) is a linear map \(f\colon K^{r} \to K^{r'}\) such that \(f(L_{fi}) \subset f(L'_{fi})\) and \(f(L_\infty) \subset L'_\infty\).\label{def:LatticePair}
	\end{enumerate}
\end{definition}

\begin{theorem} \label{thm:VBCategories}
	The four categories introduced in the definition above are equivalent.
\end{theorem}

\begin{proof}
	We describe fully faithful essentially surjective functors between the categories:
	\begin{itemize}
		\item \ref{def:OxLattice} \(\to\) \ref{def:VectorBundle}: The forgetful functor from the category of \(\Oc_X\)-lattices to that of vector bundles is clearly faithful. It is essentially surjective because any vector bundle \(E\) is isomorphic to a lattice once one fixes a basis of its generic stalk \(E_\eta\): this yields an isomorphim \(E_\eta \simeq K^r\) and we get injective maps \(\Gamma(U,E) \to K^r\) compatible with restriction maps. This yields an injective homomorphism \(E \to K_X^r\). Since a map between vector bundles induces a map between generic stalks, it is clear that this functor is full.
        \item \ref{def:OxLattice} \(\to\) \ref{def:ARLattice}: The functor sending an \(\Oc_X\)-lattice \(\mathcal{L}\) to \(\prod_{P \in M} \mathcal{L}_P\) is clearly fully faithful from the definitions of homomorphisms and the fact that a global homomorphism of \(K_X^r\) is the same thing as a linear endomorphism of \(K^r\). By the local-global principle for \(\Oc_X\)-lattices \cite[Exercise 9.16]{voight2021quaternion}, all but finitely many of the \(\mathcal{L}_P\) are equal to \(A_P\). Furthermore, all products of \(A_P\)-lattices with this property are reached. Such a product of lattices is the same thing as an \(A_R\)-lattice by \Cref{lemma:GLnR} and (essential) surjectivity follows.
        \item \ref{def:ARLattice} \(\to\) \ref{def:LatticePair}: By the local-global principle for lattices on a Dedekind domain \cite[Theorem 9.4.9]{voight2021quaternion}, if \(L = (L_P)\) is an \(A_R\)-lattice, the restriction \((L_P)_{P \in M_{fi}}\) determines a unique \(A_{fi}\)-lattice which we denote by \(L_{fi}\), furthermore every \(A_{fi}\)-lattice can bet obtained this way. We likewise define \(L_\infty\) as the unique \(A_\infty\)-lattice defined by \((L_P)_{P \in X_\infty}\). It is straightforward to deduce that sending an \(A_R\)-lattice \(L\) to the pair \((L_{fi},L_\infty)\) yields an equivalence of categories.
	\end{itemize}
\end{proof}

\begin{remark}
	By our definitions, the categories \ref{def:OxLattice},\ref{def:ARLattice} and \ref{def:LatticePair} are small (their objects form a set) and the equivalence of categories we described are not merely surjective, but they induce a bijection between the sets of objects. As a result, we may unambiguously fix an object in one of these categories and talk about the associated objects in the other two categories. If \(L\) is any type of lattice, we write \(L_X\) for the corresponding \(\Oc_X\)-lattice, \(L_R\) for the corresponding \(A_R\)-lattice and \(L_\LP\) for the corresponding lattice pair.
\end{remark}

\begin{remark}\label{remark:GammaAsCap}
    If \(L\) is an \(\Oc_X\)-lattice, we observe that \(\Gamma(U,L)\) is the subset \(\bigcap_{P \in U} L_P\) of \(\Gamma(U,K_X) = K\). Indeed, a section \(s \in \Gamma(U,L)\) has its stalk in \(L_P\) for all \(P \in U\) and is therefore sent there by the restriction maps of the sheaf \(K_X\). But these maps are the identity. Converserly, an element \(s \in \bigcap_{P \in U} L_P\) directly glue back into an element of \(\Gamma(U,L)\).
\end{remark}

\subsection{\(\Oc_X\)-lattices and \(A_R\)-lattices}

	\begin{prop}\label{thm:DoubleQuotient}
        There is a bijection between the set of isomorphism classes of rank \(r\) vector bundles and the double quotient \[GL_r(K) \backslash GL_r(R) / GL_r(A_R).\]
    \end{prop}

    \begin{proof}
	    This is an easier version of \cite[Proposition 22]{weng2018codes}. We prove the result for isomorphism classes of rank \(r\). For any \(g \in GL_r(R)\), there is an \(A_R\)-lattice \(R(g) \coloneq g(A_R^r)\). This lattice is determined by \(g\) up to an automorphism of \(A_R^r\). So, the set of \(A_R\)-lattices is in bijection with \(GL_n(R)/GL_n(A_R)\). Furthermore, two lattices in this set are isomorphic if one is the image of the other by an automorphism of \(K^r\) applied pointwise. That is, the set of isomorphism classes of \(A_R\)-lattices, and therefore of vector bundles over \(X\), is in bijection with the double quotient \(GL_r(K)\backslash GL_r(R)/GL_r(A_R)\).
    \end{proof}

    As we represent an \(A_R\)-lattice \(L\) by a matrix \(g\) such that \(L = R(g)\), we establish how properties of \(L\) may be described algebraically using \(g\). 

    \begin{definition}
        Let \(g \in GL_r(R)\), let \(L = R(g)\). We define \(\det(L) = R(\det(g))\) and \(\deg(L) = -\deg(\det(g))\).
    \end{definition}

     In order to express the tensor product of \(A_R\)-lattices as a lattice, we identify the tensor product \(R^{r} \otimes_R R^{r'}\) with \(R^{rr'}\) via the tensor product of the canonical bases. That is, if \((e_1,\hdots,e_r)\) is the canonical basis of \(R^r\) and \((e'_1,\hdots,e'_{r'})\) is that of \(R^{r'}\), we identify \(R^r \otimes R^{r'}\) with \(R^{rr'}\) via the basis \((e_1 \otimes e'_1,e_1 \otimes e'_2,\hdots,e_1 \otimes e'_{r'},e_2 \otimes e'_1,\hdots,e_r \otimes e'_{r'})\).

    \begin{prop}\label{prop:AdelicMatricesOperations}
	    Let \(g \in GL_r(R)\) and \(g' \in GL_{r'}(R)\). Let \(L = R(g)\).
        \begin{enumerate}
		\item The rank one lattice \(\det(L)\) is independent of the choice of \(g\) such that \(L = R(g)\). Furthermore, \(\det(L)_X = \det(L_X)\). 
        \item The number \(\deg(L)\) is independent of the choice of \(g\) such that \(L = R(g)\). Furthermore, \(\deg(L_X) = \deg(L)\).
	\item \(R(g) \otimes_{A_R} R(g') = R(g \otimes g')\).
	\item \(R(g) \oplus R(g') = R(g \oplus g')\).
	\item Let \(M \in M_{r',r}(K)\). Then \(M\) describes a map from \(R(g)\) to \(R(g')\) if and only if \(g'^{-1} M g \in M_{r',r}(A_R)\). That is, \(\Hom(R(g),R(g')) = M_{r',r}(K) \cap g'M_{r',r}(A_R)g^{-1}\). \label{item:AdelicHom}
        \end{enumerate}
    \end{prop}

    \begin{proof}
        The matrix \(g\) such that \(L = R(g)\) is defined up to a factor in \(GL_r(A_R)\). Such a factor has a determinant in \(A_R^\times\). Multyplying a répartition by an element of \(A_R^\times\) does not change the \(A_R\)-lattice of rank \(1\) it generates. Thus, \(\det(L)\) is independent of the choice of \(g\), and so is \(\det(L)\).

	    Then, each item is proved by observing that the result holds locally at each \(P \in M\). We note that the degree of a rank \(1\) \(A_R\)-lattice is the opposite of the degree of its generator: if \(r \in R\) is invertible and \(D = \sum_{P \in X} v_P(r_P) P\) is the associated divisor of \(r\), the line bundle associated to \(R(r)\) is in fact \(\mathcal{L}(-D)\), and by the definition of the degree of a répartition, \(\deg r = \deg D\).
    \end{proof}

    \begin{definition}
        If \(L\) is an \(A_R\)-lattice, the set
        \[L^\vee = \left\{a \in R^r \mid \forall b \in L, \sum_{i=1}^r a_i b_i \in A_R\right\}\]
        is called the dual lattice of \(L\).
    \end{definition}

    \begin{prop}\label{prop:LatticeDual}
        Let \(L\) be an \(A_R\)-lattice. The dual \(L^\vee\) is an \(A_R\)-lattice. This duality is the same as that of vector bundles:
	\[(L^\vee)_X \simeq (L_X)^\vee.\]
    \end{prop}

    \begin{proof}
	    If \(g \in GL_r(R)\) is such that \(L = R(g)\), then \(L^\vee = R({}^tg^{-1})\) and therefore \(L^\vee\) is an \(A_R\)-lattice.

        Let \(U \subset X\) be an open subset, and let \(a = (a_1,\hdots,a_r) \in \Gamma(U,L_X) \subset \Gamma(U,K_X^r) =  K^r\). For any \(s = (s_1,\hdots,s_r) \in \Gamma(U,L_X)\), define \(a(s) = \sum_{i=1}^r a_i s_i \in \Gamma(U,K_X)\). Then, for all \(P \in U\), \(a(s) \in A_P\) since \(a \in (L^\vee)_P\) and \(s \in L_P\). Thus, \(a(s) \in \Gamma(U,\Oc_X)\) seen as a subset of \(\Gamma(U,K_X)\), and \(a\) does define a homomorphism \(\Gamma(U,L) \to \Gamma(U,\Oc_X)\). It is clear from the definition of this homomorphism that it is compatible with restriction maps.

        This yields a homomorphism \(\Gamma(U,(L^\vee)_X) \to \Gamma(U,(L_X)^\vee)\). The fact that it is an isomorphism may be checked locally.
    \end{proof}

    \begin{remark} \label{rem:DualMatrix}
        We record from the proof of \Cref{prop:LatticeDual} that if \(g \in GL_r(R)\), \(R(g)^\vee = R({}^tg^{-1})\). For convenience, when \(g \in GL_r(R)\), we set \(g^\vee = {}^tg^{-1}\).
    \end{remark}
    \begin{remark} \label{rem:HomLattice}
        \Cref{item:AdelicHom} of \Cref{prop:AdelicMatricesOperations} suggests that if \(g \in GL_r(R)\) and \(g' \in GL_{r'}(R)\), then \(\Hom(R(g),R(g'))\) is the set of global sections of the \(\Oc_X\)-lattice corresponding to the free \(A_R\)-submodule \(g'M_{r',r}(A_R)g^{-1}\) of \(M_{r',r}(R)\). Any commutative ring \(B\) is self-dual when seen as a module upon itself, via the isomorphism \(b \mapsto (a \mapsto ab)\). Then, there is a natural identification \(M_{r',r}(B) \simeq (B^\vee)^r \otimes B^{r'} \simeq B^r \otimes B^{r'}\). The basis \(e_1 \otimes e'_1, e_1 \otimes e'_2, \hdots, e_2 \otimes e'_1, \hdots, e_r \otimes e'_{r'}\) which we use in general for \(B^r \otimes_B B^{r'}\) identifies with the basis of elementary matrices \((E_{11},E_{12},\hdots,E_{r1},E_{21},\hdots,E_{rr'})\).
                One checks easily that upon identifying \(M_{r',r}(R)\) with \(R^{r'r}\), \(g'M_{r',r}(A_R)g^{-1}\) is sent to \(L^\vee \otimes L'\). For this reason, we set \(\Homc(R(g),R(g')) = g'M_{r',r}(A_R)g^{-1}\) which we also identify with \(R(g)^\vee \otimes_{A_R} R(g')\).
    \end{remark}

\subsection{Cohomology of \(A_R\)-lattices}\label{sec:AdelicCohomology}
    We may give an explicit description of the cohomology of a vector bundles in terms of adeles. We first introduce some notation, and then quote a result of \cite{weng2018codes,sugahara2012adelic}, which generalises \cite[Proposition II.5.3]{serre1975groupes}.
    \begin{definition}
	    If \(L\) is an \(A_R\)-lattice or a lattice pair, we define the cohomology groups
        \[H^0(L) = L \cap K^r\]
        and
        \[H^1(L) = R^r/(L + K^r).\]
    \end{definition}

    As usual, this definition is compatible with the analogous definition on \(X\)-lattices:

    \begin{prop}\label{prop:AdelicCohomology}
        Let \(L\) be an \(A_R\)-lattice of rank \(r\). Then \(H^0(X,L_X) = \Gamma(X,L_X)\) injects in \(K^r\) and we get:
        \[H^0(L) = H^0(X,L_X).\]
        Furthermore, there is an isomorphism
        \[H^1(L) \simeq H^1(X,L_X).\]
    \end{prop}

    We first need a lemma:

    \begin{lemma}\label{lemma:ForCohomology}
	    Let \(L\) be an \(\Oc_X\)-lattice of rank \(r\). Then for any open subset \(U \subset X\), \(\Gamma(U,K_X^r/L) = \bigoplus_{P \in U} K^r/L_P\). In particular, the sheaf \(K_X/L\) is flasque.
    \end{lemma}

    \begin{proof}
	    It is enough to prove the result on an open cover of \(X\), so we may assume without loss of generality that \(U = \Spec(A)\) is an affine open subset of \(X\) over which \(L\) is free. Then, \(L_A\) may be seen as an \(A\)-lattice isomorphic to \(A^r\)r, and we must prove that \(K^r/L_A \simeq \bigoplus_{\mathfrak{p} \in \Spec(A)} K^r/L_P\). Fix \((a_1,\hdots,a_r)\) a basis of \(L_A\) and then \(K^r/L_A = \bigoplus_{i=1}^r Ka_i/Aa_i\). It is then enough to prove the result for \(L_A = A\). However, if \(A\) is a Dedekind domain and \(K\) is its fraction field, then \(K/A \simeq \bigoplus_{\mathfrak{p} \in \Spec(A)} K/A_\mathfrak{p}\) by the Chinese Remainder Theorem.
    \end{proof}

    \begin{proof}[Proof of \Cref{prop:AdelicCohomology}]
        We rephrase Serre's argument and adapt it to our more general context. By \Cref{lemma:ForCohomology}, the middle and right terms of the exact sequence
        \[0 \rightarrow L_X \rightarrow K_X^r \rightarrow K^r/L_X \rightarrow 0\]
	are flasque sheaves, so \(H^1(X,K_X^r) = H^1(X,K_X^r/L_X) = 0\) \cite[Proposition III.2.5]{hartshorne2013algebraic}. Therefore the cohomology of \(L_X\) may be computed as the kernel and cokernel of the map \(\Gamma(X,K_X^n) \rightarrow \Gamma(X,K_X^r/L_X)\). Now, \(\Gamma(K,K_X^r) = K^r\) and by \Cref{lemma:ForCohomology}, \(\Gamma(X,K_X^r/L_X) \simeq \bigoplus_{P \in X} K^r/L_P = R/L\).

    This gives isomorphisms \(H^*(L) \simeq H^*(X,L_X)\). The fact that the isomorphism of \(H^0\) groups is actually an equality under the identification of \(H^0(X,L_X)\) with its image in \(K^r\) is a consequence of \Cref{remark:GammaAsCap}.
    \end{proof}

    We conclude this subsection with an explicit description of Serre duality in our setting. This is a rephrased version of a theorem for adèles proved in \cite{sugahara2012adelic}. It would be possible to adapt Serre's proof given in \cite[Section II.8]{serre1975groupes} and prove the theorem using only the theory of répartitions and differentials. However, for efficiency purposes, we will assume that the abstract statement of Serre duality is already known and content ourselves with giving concrete formulas for computation.

    \begin{theorem}[Serre Duality]
	    Let \(\omega\) be a differential of \(K\). For any \(A_R\)-lattice \(L\) of rank \(r\), there is a perfect pairing
	    \[\begin{array}{cccc}
		    \theta_\omega \colon & H^0(\iota(\omega)^{-1} L^\vee) \times H^1(L) &\to &F \\
				 & (a,b) &\mapsto & \Res\left(\sum_{i=1}^r a_i b_i \iota(\omega)\right).
	    \end{array}\]
\end{theorem}

\begin{proof}
	We first prove that \(\theta\) is a well-defined pairing.
	For \(a \in H^0(\iota(\omega)^{-1} L^\vee)\), consider the map	
	\[\begin{array}{cccc}
		\theta'_\omega(a,\cdot) \colon & R^r &\to &F \\
				& b & \mapsto & \Res\left(\sum_{i=1}^r a_i b_i \iota(\omega)\right).
	\end{array}\]
	We prove that \(L + K^r \subset \Ker(\theta'_\omega(a,\cdot))\). Observe that \[H^0(\iota(\omega)^{-1} L^\vee) = \left\{f \in K^r \mid \forall b \in L, \sum_{i=1}^r \iota(\omega) a_i b_i \in A_R\right\}.\]
    If \(b \in L\), \(ab\iota(\omega) \in A_R^r\), so \(\theta'_\omega(a,b) = \Res\left(\sum_{i=1}^r a_i b_i \iota(\omega)\right) = 0\) since \(\Res(x) = 0\) for any \(x \in A_R\). If \(b \in K^r\), then \(\sum_{i=1}^r a_i b_i \iota(\omega) = \iota(\sum_{i=1}^r a_i b_i \omega)\) since the set of differentials of \(K\) is a \(K\)-vector space. Therefore \(\theta'_\omega(a,b) = 0\) by the Residue Theorem. It follows that \(\theta'_\omega(a,\cdot)\) factors into a unique map \(\theta_\omega(a,\cdot)\) from \(H^1(L)\) to \(k\). The pairing \(\theta_\omega\) is well defined.

	We prove that the map \(a \mapsto \theta_\omega(a,\cdot)\) is injective. Let \(a \in H^0(\iota(\omega)^{-1} L^\vee)\) Assume that \(a_{i,P} \neq 0\) for some \(i \in [r]\), \(P \in X\) and set \(n = v_P(a \iota(\omega)) + 1\), \(b_j = 0\) for \(j \neq i\), \(b_{i,Q} = 0\) for \(Q \neq P\) and \(b_{i,P} = 1/\pi_P^n\). Then \(\theta_a(b) = \Res_P(\iota(\omega) a_i b_i) \neq 0\). Therefore, the map \(\theta_a\) is non-zero over \(R^r\) and therefore over \(H^1(L) = R^r/(L + K^r)\).

	Since \(H^0(\omega^{-1} L^\vee)\) and \(H^1(L)\) are finite dimensional \(F\)-vector spaces of equal dimensions (by the abstract statement of Serre duality), it follows that the map \(a \mapsto \theta_a\) is an isomorphism. That is, \(\theta\) is a perfect pairing.
\end{proof}

\begin{remark}\label{remark:SerreNatural}
	The pairing \(\theta_\omega\) behaves naturally with the change of differential. More precisely, if \(\omega' = f\omega\) is a different differential of the field \(K\), then multiplication by \(f\) gives an isomorphism \(H^0(\iota(\omega')^{-1}L^\vee) \simeq H^0(\iota(\omega)L^\vee)\). We easily check that for any \(a \in H^0(\iota(\omega')^{-1}L^\vee)\), \(\theta_{\omega'}(a) = \theta_{\omega}(fa)\).
\end{remark}

\subsection{Extensions of \(A_R\)-lattices}

We briefly recall the general theory of extension of vector bundles. We then give an explicit construction of an extension of \(A_R\)-lattices. A reference for extensions of vector bundles is \cite[Section 7.3]{lepotier1997lectures}.

\begin{definition}
	Let \(F,G\) be vector bundles over \(X\). Then an extension of \(F\) by \(G\) is an exact sequence
	\[0 \to G \to E \to F \to 0.\]
\end{definition}
	A map of extensions is a map of exact sequences. We note that two extensions of \(F\) by \(G\) may not be isomorphic as extensions even though the vector bundles in the middle of the sequences are.
	It is well known that extensions of modules are in general classified by the cohomology group \(\Ext^1(F,G)\). In the case of vector bundles, this group is naturally isomorphic to \(H^1(X,F^\vee \otimes G)\), with the isomorphism given as follows:
    \begin{prop}\label{prop:AbstractExt}
		Let \(F\) and \(G\) be vector bundles over \(X\). Then, there is a bijection \(\delta\) between the set of isomorphy classes of extensions of \(F\) by \(G\) and \(H^1(X,F^\vee \otimes G)\). The map \(\delta\) is defined as follows: let \(\xi\) be an extension given by the exact sequence
		\[0 \to G \to E \to F \to 0.\]
		Then, the following sequence is also exact.
		\[0 \to F^\vee \otimes G \to F^\vee \otimes E \to F^\vee \otimes F \to 0.\]
		This sequence yields a map \(\partial\colon \Hom(F,F) = H^0(F^\vee \otimes F) \to H^1(F^\vee \otimes G)\). Then, \(\delta(\xi) = \partial(Id_F)\).
	\end{prop}

	This result is usually proved using injective resolutions of sheaves, which yields a construction of the map \(\delta^{-1}\). However, this is impractical in a computational setting. Instead, we adapt the methods from \cite[Section 2]{weng2018adelic} and redo the computation using the exact sequence from the proof of \Cref{prop:AdelicCohomology}.
	\begin{theorem}\label{thm:ExplicitExt}
		Let \(g' \in GL_{r'}(R)\) and \(g'' \in GL_{r''}(R)\). Let \(\kappa \in M_{r',r''}(R)\). Then \(\kappa\) represents an element of \(H^1(\Homc(R(g''),R(g')))\) and therefore an extension of \(R(g'')\) by \(R(g')\). This is represented by the exact sequence
		\[0 \to R(g') \xrightarrow{\iota} R(g) \xrightarrow{\pi} R(g'') \to 0\]
		where
		\[g = \begin{pmatrix} g' & -\kappa g'' \\ 0 & g'' \end{pmatrix} \in GL_r(R),\]
	and \(\iota\) and \(\pi\) are respectively given by the injection of \(R^{r'}\) into the \(r'\) first summands of \(R^r\) and by the projection of \(R^r\) onto its last \(r''\) summands.
	\end{theorem}

	\begin{proof}
        Recall from the proof of \Cref{prop:AdelicCohomology} that for any \(\Oc_X\)-lattice \(L\), we have the exact sequence
        \[0 \to L_X \to K^r_X \to K^r_X/L_X \to 0,\]
        which gives rise to the following long exact sequence:
        \begin{equation} H^0(L) \to K^r \to R^r/L \to H^1(L).\label{eq:LongExactSequence} \end{equation}
        Indeed, \(R^r/L\) is none other than \(\bigoplus_{P \in X} K^r/L_P\), which is the group \(H^0(X,K_X/L_X)\) by \Cref{lemma:ForCohomology}.

        Let \(L = R(g)\), \(L' = R(g')\) and \(L'' = R(g'')\).
        Writing (\ref{eq:LongExactSequence}) vertically for each term of the short exact sequence
        \[0 \to \Homc(L_X'',L'_X) \to \Homc(L''_X,L_X) \to \Endc(L''X) \to 0,\]
        we get the diagram

        \begin{tikzpicture}[baseline= (a).base]
            \node[scale=.8] (a) at (0,0){
                \begin{tikzcd}
                    0 \arrow[r] & \Hom(L'',L') \arrow[r] \arrow[d] & \Hom(L'',L) \arrow[r] \arrow[d] & \End(L'') \arrow[d] \arrow[r,"\partial"] &\hdots \\
                    0 \arrow[r] & M_{r',r''}(K) \arrow[r] \arrow[d] & M_{r,r''}(K) \arrow[r] \arrow[d] & M_{r''}(K) \arrow[r] \arrow[d] & 0 \\
                    0 \arrow[r] & M_{r',r''}(R)/\Homc(L'',L') \arrow[r] \arrow[d] & M_{r,r''}(R)/\Homc(L'',L) \arrow[r] \arrow[d] & M_{r''}(R)/\Endc(L'') \arrow[d] \\
                    \hdots \arrow[r,"\partial"] & H^1(\Homc(L'',L')) \arrow[r] & H^1(\Homc(L'',L)) \arrow[r] & H^1(\Endc(L'')) \arrow[r] & 0.
                \end{tikzcd}
            };
        \end{tikzpicture}

		At each line, the maps are between rings of matrices with coefficients either in \(K\) or in \(R\). Either way, the first map always sends a matrix \(M'\) of size \(r' \times r''\) to the matrix \(\begin{pmatrix} M \\ 0 \end{pmatrix}\) of size \(r'+r'' \times r''\) and the second map sends a matrix \(M = \begin{pmatrix} M_1 \\ M_{2} \end{pmatrix}\) of size \(r'+r'',r''\) to the matrix \(M_2\) of size \(r'' \times r''\). We respectively denote this injection and projection by \(\iota'\) and \(\pi'\), regardless of the coefficient ring.
		We wish to compute \(\delta(\xi) = \partial(Id_L'')\). By the usual proof of the snake lemma, \(\partial(Id_L'') \in H^1((L'')^\vee_R)\) is represented by a matrix \(c \in M_{r',r''}(R)\) such that there exist \(U \in M_{r',r''}(K)\) and \(V \in \Homc(L'',L) = gM_{r,r''}(A_R)g''^{-1}\) such that \(\iota'(c + U) = \begin{pmatrix} 0 \\ I_{r''} \end{pmatrix} + V\). Now, if \(M = \begin{pmatrix} M_1 \\  M_2 \end{pmatrix} \in M_{r,r''}(A_R)\) (with \(M_1\) having \(r'\) lines and \(M_2\) having \(r''\)), we get 
		\begin{align*}
			gMg''^{-1} &= \begin{pmatrix} g' & -\kappa g'' \\ 0 & g'' \end{pmatrix} \begin{pmatrix} M_1 g''^{-1} \\ M_2 g''^{-1} \end{pmatrix} \\
				   &= \begin{pmatrix} g' M_1 g''^{-1} + -\kappa g'' M_2 g''^{-1} \\ g'' M_2 g''^{-1} \end{pmatrix}.
			\end{align*}
			Therefore, setting \(M_2 = -I_{r''}\) and \(M_1 = 0\), we construct \(V = \begin{pmatrix} \kappa \\ -I_{r''} \end{pmatrix} \in \Homc(L'',L)\), and observe that \(\iota'(\kappa) = \begin{pmatrix} 0 \\ I_{r''} \end{pmatrix} + V\). It follows directly that the class \(\delta(\xi) = \partial(Id_{L''})\) in \(H^1(\Homc(L'',L'))\) is represented by the matrix \(\kappa\).
	\end{proof}

    \subsection{Restriction and conorm of an \(A_R\)-lattice}\label{sec:RestrictionConormVB}
    In this section we consider a finite separable function field extension \(K'/K\). We also set \(n = [K' : K]\). As this separable extension corresponds to a separable morphism of algebraic curves, we define the restriction and conorm of an \(A_R\)-lattice as the counterpart of respectively the direct and inverse image of a vector bundle.

    The first convention we adopt for the rest of this section is that we assume a fixed \(K\)-basis of \(K'\) denoted by \(c_1,\hdots,c_n\). Using this, we identify \(K'\) with \(K^n\) and more generally \((K')^r\) with \(K^{rn}\) as \(K\)-vector spaces. That is, if \(e_1,\hdots,e_r\) is a basis of \(L^r\), the corresponding basis of \(K^{nr}\) is \((e_1c_1,e_1c_2,\hdots,e_1c_n,e_2c_1,\hdots,e_rc_n)\). We write \(R'\) for the répartition ring \(R_{K'}\), and \(M'\) for the set of places of \(K'\). If \(Q \in M'\), we let \(Q_K\) be the place in \(M\) lying below \(Q\).

	\begin{definition} Let \(L\) be an \(A_{R}\)-lattice of rank \(r\) and let \(L'\) be an \(A_{R'}\)-lattice of rank \(r'\).
		\begin{itemize}
            \item The \emph{restriction} of \(L'\) to \(K\) is the \(A_R\)-lattice of rank \(nr'\) defined locally at \(P \in M\) by
                \[(\restriction(L'))_P \coloneqq \bigcap_{Q \mid P} L'_Q.\]
            \item The \emph{conorm} of \(L\) over \(K'\) is the \(A_{R'}\)-lattice of rank \(r\) defined locally at \(Q \in M'\) by
                \[(\conorm(L))_Q \coloneqq A_Q L_{Q_K}.\]
		\end{itemize}
	\end{definition}

    \begin{prop}
        Let \(f\colon X \to X'\) be a morphism of curves corresponding to the function field extension \(K'/K\). Let \(L\) be an \(A_R\)-lattice and \(L'\) be an \(A_{R'}\)-lattice. Then,
        \[\restriction(L')_{X'} = f_*(L'_{X'})\]
        and
        \[\conorm(L)_X = f^*(L_X).\]
    \end{prop}

    \begin{proof}
        This is directly checked on stalks.
    \end{proof}

    It is well known that for quasi-coherent sheaves, and therefore for vector bundles, there exist natural isomorphisms \(H^*(X',L'_{X'}) \simeq H^*(X,f_*(L'_{X'}))\). In our setting, this isomorphism becomes equality for \(H^0\).

Regarding \(H^1\), we will identify the space \(R^{rn}\) with its image in \(R'^r\) by the injective map \(\varphi\) defined as follows: First identify \(R^{rn}\) with the restricted product \(\left(\tilde{\prod}_{P \in M} K'\right)^r\), and then send a vector \((v_P)_{P \in M}\) to \((v_{Q})_{Q \in M'}\). We then get the following identification.

    \begin{prop}\label{prop:CohomologyEquality}
        Let \(L'\) be an \(A_{R'}\)-lattice. Then under the usual identification \((K')^r = K^{rn}\),
        \[H^0(L') = H^0(\restriction(L')).\]
        Furthermore, the map \(\varphi\) described above factors into an isomorphism \[\tilde\varphi\colon H^1(\restriction(L)) \simeq H^1(L).\]
	\end{prop}

	\begin{proof}
		The first result is proved by the direct computation:
        \begin{align*}
            H^0(L') &= H^0(X',L_{X'}) \\
                    &= \bigcap_{Q \in M'} L'_Q \\
                    &= \bigcap_{P \in M} \bigcap_{\substack{Q \in M' \\ Q \mid P}} L'_Q \\
                    &= \bigcap_{P \in M} \restriction(L')_P \\
                    &= \bigcap_{P \in M} H^0(X,\restriction(L')_X) \\
                    &= H^0(\restriction(L')).
        \end{align*}

        For the second result, we prove that \(\varphi^{-1}(L + K^{rn}) = \restriction(L) + K^r\). First, observe that \(\varphi(R^{rn})\) is the space of répartition vectors \(v\) such that \(v_Q = v_{Q'}\) if \(Q\) and \(Q'\) lie above the same place of \(K\). It is also clear that \(\varphi(K^{rn}) = (K')^r\).

		Next, we observe that \(\varphi(\restriction(L')) = L' \cap \varphi(R^{rn})\). Indeed, let \(v \in L \cap \varphi(R^{rn})\), and fix \(P \in M\). Then, for all \(1 \leq i \leq r\), the \(v_{i,Q}\) are equal for all \(Q \mid P\), and we denote their common value by \(v_{i,P}\). It follows that \[(v_{1,P},\hdots,r_{v,P}) \in \cap_{Q \mid P} L'_Q = (\restriction(L'))_P.\] Therefore, \(v \in \varphi(\restriction(L'))\) and \(L' \cap \varphi(R^{rn}) \subset \restriction(L')\). The converse inclusion is clear enough.

		As \(\varphi(K^{rn}) = (K')^r\) and \(\varphi(\restriction(L')) = L' \cap \varphi(R^{rn})\), we have \[\restriction(L') + K^{rn} = \varphi^{-1}(L' + K'^r).\] Indeed, let \(r \in R^{rn}\) such that \(\varphi(r) = s + t\), with \(s \in L'\) and \(t \in K'^r\). Set \(u \in K^{rn}\) such that \(\varphi(u) = t\) and observe that \(\varphi(r-u) \in L' \cap \varphi(R^{rn})\). So, \(r-u \in \restriction(L')\) and \(r \in \restriction(L') + K^{rn}\).

		This shows that \(\varphi\) factors into an injective map from \(H^1(\restriction(L'))\) to \(H^1(L')\). Surjectivity follows from equality of dimensions, as it is known that these two finite dimensional \(k\)-vector spaces are isomorphic from general results on quasi-coherent sheaves.
	\end{proof}

	\subsection{Indecomposable \(A_R\)-lattices}\label{sec:IndecomposableBundles}
	Since the Krull-Schmidt theorem applies to the category of vector bundles over \(X\) \cite{atiyah1956krull}, and as the direct sum of two \(A_R\)-lattices is easily characterized, we are mostly concerned with constructing vector bundles that do not split into a direct sum of vector bundles. We recall here results from \cite{arason1992indecomposable} and interpret them in terms of \(A_R\)-lattices. Results stated without proofs in this section are simple restatements of results from the sources above.

	\begin{definition}
		An \(A_R\)-lattice \(L\) is indecomposable if for any \(A_R\)-lattices \(L'\) and \(L''\) such that \(L \simeq L' \oplus L''\), either \(L'\) or \(L''\) is the zero module.

        An \(A_R\)-lattice \(L\) is absolutely indecomposable if its conorm over \(\overline{k}K\) is an indecomposable \(A_{R_{\overline{k}}}\)-lattice.
	\end{definition}

    \begin{remark}
        Since the objects of the cateogry of \(A_R\)-lattices are free \(A_R\)-modules this notion may seem trivial. However, since we restrict the maps to homomorphisms that are globally defined (that is, defined by a matrix with coefficients in \(K\)), our notion of direct sum is also restricted and there may exist indecomposable \(A_r\)-lattices of rank larger that \(1\).
    \end{remark}

	\begin{prop}[Krull-Schmidt-Atiyah]
		Any \(A_R\)-lattice \(L\) admits a decomposition into a direct sum of indecomposable \(A_R\)-lattices
        \[L \simeq \bigoplus_{i=1}^s L_i^{n_i}.\]
		Furthermore, such a decomposition is unique up to a renumbering of the summands.
	\end{prop}

	Let \(L\) be an \(A_R\)-lattice \(L\). Then, \(\End(L)\) is a \(k\)-algebra. For a \(k\)-algebra \(A\), we denote the Jacobson radical of \(A\) by \(J(A)\) and by \(D(A)\) a Wedderburn-Malcev complement of \(J(A)\) in \(A\). That is, \(A \simeq D(A) \oplus J(A)\) and \(D(A)\) is semi-simple. The complement \(D(A)\) is unique up to inner automorphisms of \(A\). We also denote by \(D(L)\) the semi-simple algebra \(D(\End(L))\). We get the following description of the structure of \(L\): 
    \begin{prop}\label{prop:KrullSchmidtEnd}
        Let \(L\) be an \(A_R\)-lattice, and let
        \[D(L) \simeq \bigoplus_{i=1}^s M_{n_i}(D_i)\]
        be the splitting of \(D(L)\) into a direct sum of simple \(k\)-algebras, where the \(D_i\) are division algebras. Then it is well known that
        \[L \simeq \bigoplus_{i=1}^s L_i^{n_i}\]
        where \(L_i\) is an indecomposable \(A_R\)-lattice and \(D(L_i) \simeq D_i\). Furthermore, the action of \(D(L)\) on \(L\) is compatible with this isomorphism. In particular, \(L\) is indecomposable if and only if \(D(L)\) is a division algebra.
	\end{prop}

	Following \cite{arason1992indecomposable}, and since the field \(k\) is perfect, we also have
    
    \begin{prop}\label{prop:AbsolutelyIndecomposable}
        An \(A_R\) lattice is absolutely indecomposable if and only if \(D(L) \simeq k\).
    \end{prop}

	In order to represent indecomposable vector bundles in terms of absolutely indecomposable vector bundles, the authors introduce the notion of trace of a vector bundle:
	\begin{definition}
		Let \(k'\) be a finite extension of \(k\), let \(X_{k'} = X \times_k \Spec(k')\) and let \(p\colon X_{k'} \to X\) be the projection map. Let \(E\) be a vector bundle over \(X_{k'}\), then the \emph{trace of \(E\)} is set to be \(\Tr_{k'/k}(E) = p_*(E)\).
	\end{definition}
	They then prove the following result:
	\begin{prop}\label{thm:IndecomposableAsTrace}
		Let \(F\) be an indecomposable vector bundle on \(X\) and let \(k'\) be a maximal field contained in \(D(F)\). Then there is an absolutely indecomposable vector bundle \(E\) on \(X_{k'}\) such that \(F = \Tr_{k'/k}(E)\).
	\end{prop}

	In order to translate \Cref{thm:IndecomposableAsTrace} in terms of \(A_R\)-lattice, we only need to give an interpratation of the trace defined above, which is easily seen to be a restriction:
    \begin{definition}
        Let \(k'/k\) be a separable extension of \(k\), and let \(k'K\) be the corresponding constant field extension of \(K\). Then if \(L'\) is an \(A_{R_{k'K}}\)-lattice, we set
        \[\Tr_{k'/k}(L') = \restriction(L').\]
    \end{definition}

\section{Explicit computations with lattice pairs}\label{sec:Algorithms}
\subsection{Algorithmic representation of lattice pairs}
Recall that if \(A\) is a Dedekind domain, and \(L\) is an \(A\)-lattice of rank \(r\), then \(L\) always admits a decomposition
\begin{equation}\label{eq:PseudoBasis}
	L = \ap_1 x_1 \oplus \ap_2 x_2 \oplus \hdots \oplus \ap_r x_r
\end{equation}
where the \(\mathfrak{a}_i\) are fractional ideals of \(A\) and the \(x_i\) lie in \(\Frac(A)^r\) \cite[Section 9.3]{voight2021quaternion}.

\begin{definition}
	A \emph{pseudo-basis} of \(L\) is a pair \(\left((\ap_1,\hdots,\ap_r),(x_1,\hdots,x_r)\right)\) satisfying \Cref{eq:PseudoBasis}.
\end{definition}

That is, \(L\) is the image of the pseudo-matrix \(\left((\ap_1,\hdots,\ap_r),\begin{pmatrix} x_1 & \hdots & x_r\end{pmatrix}\right)\) as discussed in \Cref{sec:Hermite}. We move back to our usual setting with a function field \(K\) and \(A_{fi}\) and \(A_\infty\) the finite and infinite maximal orders of \(K\).

\begin{definition}
	A matrix pair of rank \(r\) is a tuple \(g = (\ap,g_{fi},g_\infty)\), where \((\ap,g_{fi})\) is a square pseudo-matrix of size \(r\) over \(A_{fi}\) and \(g_\infty \in GL_r(K)\).
	Given such a matrix pair \(g\), we define the lattice pair \(\LP(g)\) as the pair \((PM(A_{fi}^r),g_\infty(A_\infty^r))\), and say that \(g\) is a matrix pair for \(L\) if \(L = \LP(g)\).
\end{definition}

Our purpose in this section is then to translate the results from \Cref{sec:VectorBundles} in terms of matrix representations of lattice pairs. While there is no one-to-one translation from répartition matrices to matrix pairs, a matrix pair may be represented as a répartition matrix in the following manner:
\begin{definition}
	Let \(g = (\ap,g_{fi},g_\infty)\) be a matrix pair. We call \(x_1,\hdots,x_r\) the columns of \(g_{fi}\) For any \(P \in M\), we set:
	\[
		g_P = \begin{cases}
			g_\infty \text{ if } P \in M^\infty\\
			\begin{pmatrix} \pi_P^{v_P(\ap_1)} x_1 & \hdots & \pi_P^{v_P(\ap_r)} x_r \end{pmatrix} \text{ otherwise}
		\end{cases}
	\]
	and we let \(\rep(g) = (g_P)_{P \in X}\) be the répartition matrix associated to \(g\).
\end{definition}

\begin{remark} \label{rem:AlgosEasy}
	It is clear that the \(A_R\)-lattice \(R(\rep(g))\) corresponds to the lattice pair \(\LP(g)\). Therefore, in order to give an algorithm to realize any construction discussed in \Cref{sec:VectorBundles}, it is enough to give an algorithmic construction of a matrix pair \(g\) such that \(\rep(g)\) corresponds to the same construction in terms of répartitionl matrices. In particular, any construction that is directly compatible with localizations is compatible with this correspondence.
\end{remark}

We get a first batch of straightforward constructions:
\begin{definition}\label{def:AlgosBasics}
	Let \(g = (\ap,g_{fi},g_\infty)\) and \(g' = (\ap',g'_{fi},g'_\infty)\) be matrix pairs of respective ranks \(r\) and \(r'\).
    \begin{enumerate}
        \item We define \(\det(g) = \left(\prod_{i=1}^r \ap_i,\det(g_{fi}),\det(g_\infty)\right)\).
        \item If \(I\) is a fractional ideal of \(A_{fi}\), we set \(\deg(I) = \sum_{P \in X_{fi}} v_P(I)\). If \(a \in K^\times\), we set \(\deg_{fi}(a) = \deg(aA_{fi})\) and \(\deg_\infty(a) = - \deg_{fi}(a) = \sum_{P \in X_\infty} v_P(a)\).
        \item If \(r = 1\), we set \(\deg(g) = \deg(\ap) + \deg_{fi}(g_{fi} + \deg_\infty(g_\infty)\). If \(r > 1\), we set \(\deg(g) = \deg(\det(g))\). Then, we define \(\deg(\LP(g)) = -\deg(g)\).
		\item We define \(g \otimes g' = \left((\ap_1 \ap'_1,\ap_1 \ap'_2 \hdots,\ap_1 \ap'_{r'},\ap_2\ap'_1,\hdots,\ap_r \ap'_{r'}),g_{fi} \otimes g'_{fi},g_\infty \otimes g'_\infty \right)\).
        \item We define \(g \oplus g' = \left((\ap_1,\hdots,\ap_r,\ap'_1,\hdots,\ap'_{r'}), g_{fi} \oplus g'_{fi}, g_\infty \oplus g'_\infty \right)\).
        \item We define \(g^\vee = \left((\ap_1^{-1},\hdots,\ap_r^{-1}), (g_{fi}^t)^{-1},(g_\infty^t)^{-1}\right)\).
    \end{enumerate}
\end{definition}

\begin{theorem}\label{thm:BasicAlgos}
	Let \(g = (\ap,g_{fi},g_\infty)\) and \(g' = (\ap',g'_{fi},g'_\infty)\) be matrix pairs of respective ranks \(r\) and \(r'\).
	\begin{enumerate}
		\item \(\rep(\det(g)) = \det(\rep(g))\). \label{algo:Determinant}
        \item \(\deg(g) = \deg(\rep(g))\).
		\item \(\rep(g \otimes g') = \rep(g) \otimes \rep(g')\). \label{algo:TensorProduct}
		\item \(\rep(g \oplus g') = \rep(g) \oplus \rep(g')\).\label{algo:DirectSum}
		\item \(\rep(g^\vee) = \rep(g)^\vee\).\label{algo:Dual}
	\end{enumerate}
\end{theorem}

\begin{proof}
	All of these construction may be checked locally. One must simply check that the operation done on the tuple of ideals match the movements of the columns of \(g_{fi}\).
\end{proof}

\begin{remark}\label{rem:AlgoHom}
	It is a bit more tedious to directly translate our statement on homomorphisms of lattices. Instead, we may simply define \(\Homc(L,L') = L^\vee \otimes L'\) and recall the isomorphism \(M_{r',r}(K) \simeq K^{rr'}\) given by the basis of elementary matrices. Then, an algorithm for computing the lattice pair of homomorphisms follows from \Cref{thm:BasicAlgos}.
\end{remark}

\begin{example}\label{ex:LBundle}
    We let \(k = \F_{7}\) and consider the genus \(1\) function field \(K = k(x,y)/(y^2 - x^3 - x)\). We let \(\pi = \frac{y}{x^2}\) be a local uniformizer at infinity. Observe that \(\mathfrak{p} = \langle x,y\rangle\) is a prime ideal of \(A_{fi}\). We consider the lattice pair
    \[L = \LP\left(\left(A_{fi},\mathfrak{p}^{-1}\right),\begin{pmatrix} \frac{x^2}{x^2 + 4} & 0 \\ 0 & 1\end{pmatrix}, \begin{pmatrix} 1 & -\pi^{-1} \\ 0 & 1 \end{pmatrix}\right).\]
    We compute
    \[\det(L) = \LP\left(\mathfrak{p}^{-1},\frac{x^2}{x^2 + 4},1\right),\]
    and therefore
    \[\deg(L) = -\deg(\mathfrak{p}^{-1}) = 1.\]
\end{example}

\subsection{Restriction and conorm of a lattice pair}\label{sec:DirectImages}
We adopt the same notations and setting as in \Cref{sec:RestrictionConormVB}.

\begin{definition}
    Let \(L' = (L'_{fi},L'_\infty)\) be a lattice pair of rank \(r\) over \(K'\). We define \(\restriction(L')\) as the pair \((\restriction(L'_{fi}),\restriction(L'_\infty))\), where \(\restriction(L'_*)\) is the lattice \(L'_*\) seen as an \(A_*\)-lattice under the identification \(K'^r = K^{rn}\) (where \(*\) is either \(fi\) or \(\infty\)).

    Let \(L = (L_{fi},L_\infty)\) be a lattice pair of rank \(r\) over \(K\). We define \(\conorm(L)\) as the pair \((\conorm(L_{fi}),\conorm(L_\infty))\), where \(\conorm(L_*) = A'_* L_* \subset K'^r\) is an \(A'_*\)-lattice.
\end{definition}

One checks readily that these definitions are compatible with the equivalent definitions on \(A_R\)-lattices.

Matrix pairs for \(\restriction(L)\) and \(\conorm(L)\) may easily be computed.
\begin{definition}
    Let \(g = (\ap,g_{fi},g_\infty)\) be a matrix pair defined over \(K\). We set \[\conorm(g) = \left((\ap_1 A'_{fi},\hdots,\ap_r A'_{fi}),g_{fi},g_\infty\right).\] 
\end{definition}

The definition of the restriction of a matrix-pair is a bit more tedious to write down. We first give the definition of the restriction of a pseudo-matrix. The definition is given over any Dedekind domain with fraction field \(K'\), as it applies to both \(A_{fi}\) and \(A_\infty\), with the specificity that we only consider pseudo matrices with trivial coefficient ideals over \(A_\infty\), since it is a PID.

\begin{definition}
    Let \(A'_*\) be a Dedekind domain with fraction field \(K'\), and set \(A_* = A'_* \cap K\). Let \(PM = (\ap,g)\) be a pseudo-matrix of rank \(r\) over \(A'_*\). The ideals \(\ap_i\) each admit a pseudo-basis \((\bp_{i1},\hdots,\bp_{in}),(a_{i1},\hdots,a_{in})\) over \(A_*\).

    Then, we define the pseudo-matrix \(\restriction(PM)\) of rank \(rn\) over \(A_*\) with coefficient ideals \((\bp_{11},\bp_{12},\hdots,\bp_{1n},\bp_{21},\hdots,\bp_{rn})\) and matrix
    \[\begin{pmatrix} a_{11} g_{11} & a_{12} g_{11} & \hdots & a_{1n} g_{11} & a_{21} g_{12} & \hdots & a_{rn} g_{1r} \\
        a_{11} g_{21} & a_{12} g_{21} & \hdots & a_{1n} g_{21} & a_{21} g_{22} & \hdots & a_{rn} g_{2r} \\
        \vdots & \vdots & \ddots & \vdots & \vdots & \ddots & \vdots \\
        a_{11} g_{r1} & a_{12} g_{r1} & \hdots & a_{1n} g_{r1} & a_{21} g_{r2} & \hdots & a_{rn} g_{rr}
    \end{pmatrix},\]
    where each \(a_{ij} g_{\ell k}\) is understood as a column vector in \(K^n\) representing an element of \(K'\) in the usual fixed basis.

    Now, if \(g' = (\ap',g'_{fi},g'_\infty)\) is a matrix pair over \(K'\), we may define
    \[\restriction(g') = (\ap,g_{fi},g_\infty),\]
    where \((\ap,g_{fi}) = \restriction((\ap',g'_{fi}))\), and likewise \(g_\infty = \restriction(g'_\infty)\), where it is understood that all coefficient ideals of \(g'_\infty\) are equal to \(A'_\infty\), which admits a basis over \(A_\infty\).
\end{definition}

\begin{theorem}\label{thm:ResConorm}
        Let \(g\) be a matrix pair over \(K\). Then,
		\[\conorm(\LP(g)) = \LP(\conorm(g)).\]
		Let \(g'\) be a matrix pair over \(K'\). Then, \[\restriction(\LP(g')) = \LP(\restriction(g')).\]
\end{theorem}

\begin{proof}
	The first claim is clear. For the second one, the definition of the matrix pair \(\restriction(g')\) is simply an explicit writing of pseudo-bases of lattices \(L'_{fi}\) and \(L'_\infty\) in \(K'^r\) identified with \(K^{nr}\).
\end{proof}

This allows us to construct traces of vector bundles as defined in \Cref{sec:IndecomposableBundles}, but also to express the restriction to \(k(x)\) of a lattice pair, which will be a key tool in the computation of global sections.

\begin{example}\label{ex:Restriction}
    We compute the restriction \(\restriction(L)\) over \(k(x)\) of the lattice pair \(L\) from \Cref{ex:LBundle}. A basis of \(A_{fi}\) over \(k[x]\) is \((1, y)\), a \(k[x]\)-basis of \(\mathfrak{p}^{-1}\) is \((1, \frac{y}{x})\) and a basis of \(A_\infty\) over the valuation ring at infinity of \(k(x)\) is \((1, \pi) = (1, \frac{y}{x^2})\). It follows that \(\restriction(L) = \left(\ap,g_{fi},g_\infty\right)\), with
    \[\ap = (k[x], k[x], k[x], k[x]),\]
    \[g_{fi} = \begin{pmatrix}\frac{x^2}{x^2 + 4} & 0 & 0 & 0 \\
        0 & \frac{x^2}{x^2 + 4} & 0 & 0 \\ 
        0 & 0 & 1 & 0 \\
    0 & 0 & 0 & \frac{1}{x}\end{pmatrix},\]
and
    \[g_\infty = \begin{pmatrix}
        1 & 0 & 0 & -1 \\
        0 & \frac{1}{x^2} & \frac{-x}{x^2+1} & 0 \\
        0 & 0 & 1 & 0 \\
        0 & 0 & 0 & \frac{1}{x^2}
    \end{pmatrix}.\]
\end{example}
\subsection{Computing cohomology groups and extensions}
If \(L\) is a lattice pair, we define \(H^i(L)\) as \(H^i(L_R)\). Given a matrix pair \(g\), we aim to compute \(k\)-bases for the spaces \(H^0(\LP(g))\) and \(H^1(\LP(g))\).

\subsubsection{Computing global sections of a lattice pair}

The computation of \(H^0(L)\) relies on the following simple observation:
\begin{lemma}
	Let \(L=(L_{fi},L_\infty)\) be a lattice pair. Then \(H^0(L) = L_{fi} \cap L_\infty\).
\end{lemma}

\begin{proof}
	This is clear using \Cref{remark:GammaAsCap}.
\end{proof}

First, we assume that \(K =k(x)\). In this case, note that \(A_{fi} = k[x]\) is a PID and every projective \(A_{fi}\)-module is free. Therefore, we omit the tuple of ideals \(\ap\) in every matrix pair and simply assume that all ideals involved are equal to \(A_{fi}\). Then, the computation of the intersection \(L_{fi} \cap L_\infty\) reduces to matrix reduction as discussed in \Cref{sec:Popov}. The method discussed here is adapted from \cite[Lemma 25]{ivanyos2018computing} and \cite{hess2002computing}.

Let \((g_{fi},g_\infty)\) be a matrix pair over \(K\). Then, the matrix pair \(g' = (g_\infty^{-1} g_{fi},Ir)\) represents an isomorphic lattice pair, and furthermore \(H^0(\LP(g)) = g_\infty\left(H^0(\LP(g'))\right)\). Upon applying the global isomorphism \(g_\infty^{-1}\), we may assume without loss of generality that \(g_\infty\) is the identity matrix.

Then, a vector \(v \in K^r\) lies in \(L_\infty\) if and only it \(|v| \leq 0\) (see \Cref{def:NormVector}). Assume that \(e_1,\hdots,e_r\) is a reduced basis of \(L_{fi}\), in the sense that the matrix \(\begin{pmatrix} e_1 & \hdots & e_r \end{pmatrix}\) is reduced. Then, by \Cref{prop:DegreePredictability}, \(\sum_{i=1}^r a_i e_i \in L_\infty\) if and only if \(\deg(a_i) \leq -|e_i|\) for all \(1 \leq i \leq r\). Since \(\sum_{i=1}^r a_i e_i \in L_{fi}\) if and only if \(a_i \in k[x]\) for all \(1 \leq i \leq r\), a basis of \(H^0(\LP(g))\) is \[(x^j e_i)_{\substack{1 \leq i \leq r \\ 0 \leq j \leq -|e_i|}}.\]

We recall the algorithm as \Cref{algo:H0} for the convenience of the reader. Note that if \(r = \frac{P}{Q} \in k(X)\) such that \(P,Q \in k[x]\), we define \(\deg r\) as \(\deg P - \deg Q\).
\begin{algorithm}
	\KwIn{a matrix pair \(g = (g_{fi},g_\infty)\) over k(X)}
	\KwOut{A \(k\)-basis of \(H^0(\LP(g))\)}
	Set \(M = (g_\infty)^{-1}g_{fi}\)\;
	Compute \(d \in k[x]\) such that \(dM \in M_{rn}(k[x])\)\;
    Compute a reduced basis \(\mathcal{B} = (b_1,\hdots,b_{rn})\) of the \(k[x]\)-lattice generated by the columns of \(dM\)\;
\KwRet{\(\left\{\frac{x^j}{d} g_\infty(b_i) : 1 \leq i \leq rn \text{ and } 0 \leq j \leq \deg(d) -|b_i|\right\}\)}\;
	\caption{Computing the global sections of a lattice pair over \(\Pp_k^1\).}
	\label{algo:H0}
\end{algorithm}

\begin{theorem}\label{thm:H0}
	If \(K = k(x)\), \Cref{algo:H0} outputs a \(k\)-basis of \(H^0(\LP(g))\). If \(k\) is a finite field, then \Cref{algo:H0} runs in polynomial time.
\end{theorem}

\begin{proof}
    The correctness of \Cref{algo:H0} has already been discussed above. Since there exist efficient algorithms for computing a reduced basis (see \Cref{sec:Popov}), and since the size of the output is at most \(r (\deg(d)+1)\), the algorithm runs in polynomial time.
\end{proof}

\begin{cor}\label{cor:H0}
    For a general separable extension \(K/k(x)\) and a matrix pair \(g\), a basis of \(H^0(\LP(g))\) may be computed in polynomial time.
\end{cor}

\begin{proof}
    First, compute \(H^0(\restriction(\LP(g)))\) using \Cref{algo:H0}. Then, applying \Cref{prop:CohomologyEquality}, a basis of \(H^0(\restriction(\LP(g)))\) is a basis of \(H^0(\LP(g))\) upon the identification \(K^r = k(x)^{rn}\). A representation of the vectors of the basis in \(K^r\) may be computed using the basis \(1,y,\hdots,y^{n-1}\).
\end{proof}

\begin{remark}
    In \Cref{algo:H0} we may compute the Popov normalized form of the matrix \(M\) instead of a mere reduced equivalent matrix if we want the algorithm to output a more predictable basis of \(H^0(\LP(g))\).
\end{remark}

\begin{example}\label{ex:H0}
    We gather again notations from \Cref{ex:LBundle,ex:Restriction}. Compute
    \[M \coloneq g_\infty^{-1} g_{fi} = \begin{pmatrix}
        \frac{x^2}{x^2 + 4} & 0 & 0 & 0 \\
        0 & \frac{x^4}{x^2 + 4} & 0 & 0\\
        0 & \frac{x^3}{x^2 - 1} & 1 & 0\\
        x & 0 & 0 & x
    \end{pmatrix}.\]
    We compute the Popov form of its numerator, and obtain the reduced form
    \[M' = \begin{pmatrix}
        \frac{x^2}{x^2 + 4} & 0 & 0 & \frac{4x}{x^2 + 4}\\
        0 & \frac{x^3}{x^2 - 1} & \frac{2x^4}{x^4 + 3x^2 + 3} & 0 \\
        0 & 1 & -x & 0\\
        0 & 0 & 0 & x
    \end{pmatrix}.\]
    It follows that a basis of \(H^0(\restriction(L))\) is
    \[g_\infty \begin{pmatrix} \frac{x^2}{x^2 + 4} \\ 0 \\ 0 \\ 0 \end{pmatrix} = \begin{pmatrix} \frac{x^2}{x^2 + 4} \\ 0 \\ 0 \\ 0 \end{pmatrix}.\]
    Therefore, a basis of \(H^0(L)\) is
    \[\begin{pmatrix} \frac{x^2}{x^2 + 4} \\ 0\end{pmatrix}.\]
\end{example}

We also observe that given an element \(f \in H^0(L)\) for some lattice pair \(L\), we may compute the coordinates of \(f\) in terms of a given basis of \(H^0(L)\) (for instance the one computed by the algorithm of \Cref{cor:H0}).

\begin{lemma}\label{lemma:LinearAlgebraH0}
    Let \(L\) be a lattice pair of rank \(r\) and let \(f \in H^0(L) \subset K^r\). Let \(m_1,\hdots,m_s\) be a \(k\)-basis of \(H^0(L)\). We may compute in polynomial time a vector \(a \in k^s\) such that \(f = \sum_{i=1^s} a_i m_i\).
\end{lemma}

\begin{proof}
    Fix a place \(P\) of \(K\) and a local uniformizer \(\pi_P\) at \(P\). For each \(1 \leq i \leq r, 1 \leq j \leq s\), let \(m_{ij}\) be the \(i\)th component of the basis vector \(m_i\). We write \(v_i = \min_{1 \leq j \leq s} v_P(e_{ij})\) (if all the \(e_{ij}\) are zero, simply set \(v_i = 0\)). Then, for any \(\ell \in \N\) we define the map
    \[\begin{array}{rccc}
        \varphi_{P,\ell}\colon &K^r &\to &k^{r\ell} \\
                               &f &\mapsto & (f_i^{(v_i + j)})_{\substack{1 \leq i \leq r \\ 0 \leq j \leq \ell - 1}}.
    \end{array}\]
    Now, consider the matrix \(N_{P,\ell}\) of size \(r\ell \times s\) whose columns are the \(\varphi_{P,\ell}(m_j)\). The matrix \(N_{P,\ell}\) has rank \(s\) if and only if the restriction of \(\varphi_{P,\ell}\) to \(H^0(L)\) is injective, and the coordinates of an element \(f \in H^0(L)\) with respect to basis \(m_1,\hdots,m_s\) may be computed as a vector \(a \in k^s\) such that \(N_{P,\ell} a = \varphi_{P,\ell}(f)\).

    All that is left is to prove that the restriction of \(\varphi_{P,\ell}\) to \(H^0(L)\) is injective for some \(\ell\) bounded by a polynomial in the size of the input. Let \(f = \sum_{i=j}^s a_j m_j \in H^0(L)\). Then, if \(1 \leq i \leq r\) and \(f_i \neq 0\), \(\height(f_i) \leq \sum_{j=1}^s \height(m_{ij})\), and thus \(v_P(f_i) \leq \sum_{j=1}^s \height(m_{ij})\). It follows that if \(\ell > \max_{1 \leq i \leq r} v_i + \sum_{j=1}^s \height(m_{ij})\), the map \(\varphi_{P,\ell}\) is injective over \(H^0(L)\). Since \(v_i\) is itself bounded by \(\max_{1 \leq j \leq s} \height(m_{ij})\), \(\ell\) may indeed be chosen of polynomial size in the input.
\end{proof}

\subsubsection{Computing the group \(H^1\)}\label{sec:H1}

By Serre duality, computing the group \(H^1(L)\) for a lattice pair \(L\) can be done by computing the \(k\)-vector space \(H^0(\iota(\omega)^{-1} L^\vee)\) for some differential \(\omega\). However, for applications such as computing extensions of vector bundles, it is desirable to be able to represent elements of \(H^1(L)\) as elements of \(R^r\). For the rest of the section, we assume that the constant base field \(k\) is \emph{finite}.

Our strategy will be to use the linearization technique used in \(\Cref{lemma:LinearAlgebraH0}\) to turn the inversion of the Serre duality map into a linear equation.

Fix \(Q_0 \in M^\infty\), and a local uniformizer \(\pi_{Q_0}\) such that \(v_Q(\pi_{Q_0}) = 0\) for \(Q \in M^\infty \setminus \{Q_0\}\). Let \(k_0\) be the residue field of \(Q_0\) and let \(\omega = d(\pi_{Q_0})\). For any integer \(\ell\), we write \(\lceil \ell \rceil_0\) for the smallest power of \(|k_0|\) larger or equal to \(\ell\). That is, \(\lceil \ell \rceil_0 = |k_0|^{\lceil \log(\ell)/\log(|k_0|) \rceil}\). 

We present \Cref{algo:H1} which, given a basis of \(H^0(\iota(\omega)^{-1}L^\vee)\) of size \(s\) and a linear form represented in this basis by a row vector \(\varphi \in k^s\), outputs a vector \(v \in K^r\) such that the infinite répartition vector \(v_\infty\) satisfies \(\theta_\omega(\cdot,v_\infty) = \varphi\).

\begin{algorithm}
	\KwIn{A matrix pair \(g = (\ap,g_{fi},g_\infty)\) over \(X\)}
	\KwIn{A matrix \(M = (m_{i,j}) \in M_{r,s}(K)\) whose columns are a \(k\)-basis of \(H^0(\iota(\omega)^{-1}L^\vee)\)}
	\KwIn{A row vector \(\varphi = \begin{pmatrix} 1 & \hdots & f_s \end{pmatrix} \in M_{1,s}(k)\) representing a linear form on \(H^0(\iota(\omega)^{-1}L^\vee)\) written in the basis given by \(M\)}
	\KwOut{\(a \in K^r\) such that the linear form represented by \(f\) is \(\theta_\omega(\cdot,a_\infty)\)}
	For \(1 \leq i \leq r\) and \(Q \in X_\infty\), set \(v^Q_i = \min_{1 \leq j \leq s}\left(v_P(m_{i,j})\right)\)\;
    Compute the matrix \(N_{Q_0,\ell}\) (see \Cref{lemma:LinearAlgebraH0}) for increasing values of \(\ell\) until it has rank \(s\)\;
    Let \(R = \begin{pmatrix}a_1 & \hdots & a_{r\ell}\end{pmatrix} \in M_{1,r\ell}\left(k_0\right)\) such that \(\Tr_{k_{Q_0}/k}(RN) = \varphi\)\;
    Let \(\widetilde{R} = \begin{pmatrix} \widetilde{a}_1 & \hdots & \widetilde{a}_{r\ell} \end{pmatrix} \in M_{1,r\ell}(A_{Q_0})\) be a lift of \(R\) in \(M_{1,r}(A_{Q_0})\) such that \(v_Q(\widetilde{a}_i) \geq 1\) for all \(1 \leq i \leq r\ell\) and \(Q \in M^\infty \setminus \{Q_0\}\)\;
    Let \(\pi \in K\) such that \(\lceil \ell \rceil_0 v_Q(\pi) \geq \max_{1 \leq i \leq r} -v_i^Q\) for all \(Q \in X_\infty \setminus \{Q_0\}\) and \(v_{Q_0}(\pi-1) \geq 1\)\;
    \KwRet{\(\left(\pi^{\lceil \ell \rceil_0} \pi_{Q_0}^{-v_i^{Q_0}-1}\sum_{j=0}^{\ell-1} \pi_{Q_0}^{-j} \widetilde{a}_{(i-1)\ell+j+1}^{\lceil \ell \rceil_0} \right)_{1 \leq i \leq r}\)}
	\caption{Representing elements of \(H^1(L)\)}
	\label{algo:H1}
\end{algorithm}

\begin{theorem}\label{thm:H1}
	Assume that the constant base field \(k\) is finite. \Cref{algo:H1} is correct and terminates after performing a polynomial amount of arithmetic operations in \(k\).
\end{theorem}

\begin{proof}
    First, observe that \Cref{algo:H1} terminates in polynomial time: each line of the algorithm corresponds either to linear algebra over \(k\) or to a task discussed in \Cref{sec:BasicAlgos}. 

    We prove that the output of the algorithm is correct. Set \(c=(c_1,\hdots,c_r)\) as the coordinates of the output. If \(M_j\) is the vector given as the \(j\)-th column of \(M\), we claim that \(\theta_\omega(M_j,c_\infty) = \varphi_j\). 

    Now, \(\theta_\omega(M_j,c_\infty) = \sum_{i=1}^r \sum_{Q \in X_\infty} \res_{\pi_Q}(m_{i,j}c_i\omega)\). And the result will follow from the identity \(RN = \varphi\) if we prove that for \(1 \leq i \leq r\) and \(1 \leq j \leq s\), setting 
    \[\mu_{ij} = \left(\pi^{\lceil \ell \rceil_0} \pi_{Q_0}^{-v_i^{Q_0}-1}\left(\sum_{\alpha=0}^{\ell-1} \pi_{Q_0}^{-\alpha} \widetilde{a}_{(i-1)\ell+\alpha+1}^{\lceil \ell \rceil_0} \right)m_{i,j} \omega \right),\]
    we have
    \[\mu_{ij}^{(-1)} = \sum_{\alpha=0}^{\ell-1} a_{(i-1)\ell + \alpha + 1} m_{i,j}^{(-v_i+\alpha)}\]
    and \(v_Q(\mu_{ij}) \geq 0\) for \(Q \in M^\infty \setminus \{Q_0\}\)
    where, for any \(a \in K\) and integer \(n\), \(a^{(n)}\) is the coefficient of degree \(n\) in the expansion of \(a\) as a formal series in variable \(\pi_{Q_0}\). We fix \(1 \leq i \leq r\) and \(1 \leq j \leq s\).

    Let \(Q \in M^\infty \setminus \{Q_0\}\). By construction, \(v_Q(\pi^{\lceil \ell \rceil_0}) \geq \max(-v_i^Q)\) and it follows readily that \(v_Q(\mu_{ij}) \geq 0\) since \(v_Q(\pi_{Q_0}) = 0\) and \(\widetilde{a}_m \in A_Q\) for all \(1 \leq m \leq rn\).

    Now, we have the following:
    \[\pi^{\lceil \ell \rceil_0} = 1 + O(\pi_{Q_0})^\ell\]
    and
    \[\widetilde{a}_i^{\lceil \ell \rceil_0} = a_i + O(\pi_{Q_0}^\ell).\]

    Then, we get
    \[\mu_{i,j} = \left(\sum_{\alpha = 0}^{\ell-1} a_{(i-1)\ell + \alpha + 1} m_{i,j}^{(-v_i + \alpha)} \right) \pi_{Q_0}^{-1} + O(1).\]
\end{proof}

\begin{cor}\label{cor:ExtensionsAlgorithm}
    There exists a polynomial time algorithm which, given matrix pairs \(g'\) and \(g''\), as well as a row vector \(\varphi\) representing a \(k\)-linear form over \[H^0\left(\iota(\omega)^{-1}\Homc(g',g'')\right),\] returns the corresponding extension of \(\LP(g'')\) by \(\LP(g')\).
\end{cor}

\begin{proof}
    Let \(r',r''\) be the respective orders of \(g'\) and \(g''\). Using \Cref{algo:H1}, one may compute \(\kappa \in K^{r' \times r''}\) such that the infinite répartition matrix \(\kappa\) represents the element of \(H^1\left(\Homc(R(g''),R(g'))\right) = \Ext^1\left(R(g''),R(g')\right)\) corresponding to \(\varphi\).

    Then, adapting \Cref{thm:ExplicitExt}, the corresponding extension is given by the matrix pair \((\ap,g_{fi},g_\infty)\) with
    \[\bp = (\ap'_1,\hdots,\ap'_r,\ap''_1,\hdots,\ap''_{r'}),\]
    \[g_{fi} = \begin{pmatrix} g'_{fi} & (0) \\
                (0) & g''_{fi} \end{pmatrix},\]
    and
    \[g_\infty = \begin{pmatrix} g'_\infty & -\kappa g''_\infty \\
            (0)  & g''_\infty \end{pmatrix}.\]
\end{proof}

\begin{example}\label{ex:H1}
    Let \(L\) be again as in \Cref{ex:LBundle,ex:Restriction,ex:H0}
    Since \(\deg(L) = 1\), we get \(H^1(L) = 0\). Instead, we compute \(H^1(L^\vee)\). We let \(\omega = d\pi\). Since the field \(K\) has genus \(1\), the differential \(\omega\) has a principal divisor. It is the divisor of \(\frac{x^2 + 3}{x^2}\), and we may represent \(\iota(\omega)^{-1}\) by the following matrix pair of size \(1\):
    \[\left(A_{fi},\frac{x^2}{x^2+3},1\right).\]
    Now, \(H^1(L^\vee) = H^0(\iota(\omega)^{-1}L)\). Applying \Cref{algo:H0}, we compute a basis for the \(k\)-vector space \(H^0(\iota(\omega)^{-1}L)\). We find that it has dimension \(1\) and is generated by
    \[v \coloneq \begin{pmatrix} \frac{x^4}{x^4 + 5} \\ 0\end{pmatrix}.\]
    Since \(K\) has only one place at infinity and \(\dim_k(H^0(\iota(\omega)^{-1}L)) = 1\), applying \Cref{algo:H1} is straightorward: as \(\frac{x^4}{x^4 + 5} = 1 + O(\pi^8)\), 
    \[\Res_\infty\left(v,\begin{pmatrix}\pi^{-1} \\ 0 \end{pmatrix}\right) = 1.\]
    Then, the element of \(H^1(L^\vee)\) dual to \(v\) is represented by the infinite répartition vector \(\begin{pmatrix} \pi_\infty^{-1} \\ 0 \end{pmatrix}\).
\end{example}

\subsection{Computing isomorphisms between lattice pairs}
    We wish to be able to decide whether two lattice pairs are isomorphic and, if they are, find an isomorphism. We give a probabilistic algorithm of the Monte-Carlo type for this task when the field \(k\) is large enough, and a deterministic algorithm for a weakening of the problem (the lattice pairs are assumed indecomposable) when \(k\) is any finite field.

    \begin{theorem}\label{thm:AlgoIsomLargeField}
        Let \(L,L'\) be lattice pairs such that \(\dim_k \End(L) = \dim_k \Hom(L,L') = \dim_k \Hom(L',L)\). Let \(s\) be the dimension of these spaces and we assume that \(|k| > s\). There is a polynomial time Monte-Carlo algorithm which outputs an isomorphism \(\varphi\colon L \to L'\) with probability at least \(1 - s/|S|\), where \(S\) is a subset of \(k\) in which we can sample random elements.
    \end{theorem}

    \begin{proof}
        First, observe that we may compute \(\End(L), \End(L')\) and \(\Hom(L,L')\) by applying \Cref{cor:H0} to the lattice pairs \(L^\vee \otimes L, L^\vee \otimes L'\) and \((L')^\vee \otimes L\). Their elements are represented as matrices in \(M_r(K)\), and the matrix products gives a bilinear map from \(\Hom(L,L') \times \Hom(L',L)\) to \(\End(L)\). This, together with a fixed choice of bases of \(\Hom(L',L)\) and \(\End(L)\) gives a map \(\alpha\colon \Hom(L,L') \to M_s(k)\). Observe that \(f \in \Hom(L,L')\) is an isomorphism if and only if \(\alpha(f)\) is an invertible matrix. That is, if and only if \(\det(\alpha(f)) = 0\).

        Now, setting \(f = \sum_{i=1}^s a_i m_i\), where \((m_i)\) is a basis of \(\Hom(L,L')\), we see that \(\deg(\alpha(f))\) is a homogeneous polynomial of degree \(s\) in the \(a_i\). By the Schwartz-Zippel lemma, if \(S \subset k\) is a subset of size at least \(s+1\), the probability that a uniform random element of \(\bigoplus_{j=1}^s S m_j\) is an isomorphism is at least \(1 - s/|S|\). The algorithm is then simply to sample a random element of \(\Hom(L,L')\).
    \end{proof}

    When \(k\) is a finite field, the approach of \Cref{thm:AlgoIsomLargeField} does not work if \(s\) is too large, as the Schwarts-Zippel lemma fails. For now, we only give an algorithm for the case that \(L\) is an indecomposable lattice pair. This will be used as a subroutine in \Cref{algo:Splitting}, which will then be used to compute isomorphisms in the general case (see \Cref{cor:IsomFiniteField}).

    \begin{algorithm}
        \KwIn{Matrix pairs \(g\) and \(g'\) of rank \(r\) such that \(\LP(g)\) is indecomposable}
        \KwOut{A matrix \(T \in M_r(K)\) giving an isomorphism from \(\LP(g)\) to \(\LP(g')\) if \(\LP(g) \simeq \LP(g')\), and \(\bot\) otherwise}
        Compute structure constants for the \(k\)-algebra \(A = \End(\LP(g \oplus g'))\)\;\nllabel{line:StructureConstants}
        Compute sub-algebras \(S\) and \(R\) such that \(A = S \oplus R\), \(S\) is semi-simple, and \(R\) is the Jacobson radical of \(A\)\; \nllabel{line:SplittingEnd}
        \If{\(S\) is not simple}{ \nllabel{line:CheckSimple}
            \KwRet{\(\bot\)}
        }
        Compute \(s \in \N\), a finite extension \(k'/k\), and an isomorphism \(\varphi\colon S \simeq M_{s}(k')\)\; \nllabel{line:MatrixSplit}
        \If{\(s \neq 2\)}{
            \KwRet{\(\bot\)}\;
        }
        Compute \(P \in GL_{2s}(k')\) such that \(P\varphi(\mathrm{Id}_\LP(g))P^{-1} = \begin{pmatrix} 1 & 0 \\ 0 & 0 \end{pmatrix}\) and \(P\varphi(Id_{\LP(g')})P^{-1} = \begin{pmatrix} 0 & 0 \\ 0 & 1 \end{pmatrix}\)\; \nllabel{line:IdemDiagonalization}
        \KwRet{\(\varphi^{-1}\left(P^{-1}\begin{pmatrix} 0 & 0 \\ 1 & 0 \end{pmatrix}P\right)\)}
        \caption{Computing isomorphisms between quasi-indecomposable lattice pairs over finite fields}
        \label{algo:IsomIndecomposable}
    \end{algorithm}

    \begin{lemma}\label{lemma:AlgoIsomIndecomposable}
        Assume that \(k\) is a finite field. Then \Cref{algo:IsomIndecomposable} is correct and runs in polynomial time.
    \end{lemma}

    \begin{proof}
        We first discuss the algorithm line by line, proving that the task may be done in polynomial time.
        \begin{description}
            \item[\Cref{line:StructureConstants}:] A basis of \(\End(\LP(g \oplus g'))\) may be computed using \Cref{cor:H0}. Then, the structure constants may be computed using \Cref{lemma:LinearAlgebraH0}.
            \item[\Cref{line:SplittingEnd}:] The Jacobson radical of \(A\) may be computed using \cite[Theorem 2.7]{ronyai1990computing} and a basis of \(S\) may be computed using \cite[Theorem 3.1]{degraaf1997computing}.
            \item[\Cref{line:CheckSimple}:] Checking that \(S\) is simple can be done by checking if the center of \(S\) is a field. See \cite[Section 3]{ronyai1990computing}.
            \item[\Cref{line:MatrixSplit}:] This may be done using \cite[Theorem 5.2]{ronyai1990computing}. Note that since a finite field has trivial Brauer group, a simple \(k\)-algebra is always of the form \(M_{n_i}(k')\), where \(k'\) is a finite extension of \(k\).
            \item[\Cref{line:IdemDiagonalization}:] Observe that the matrices \(\varphi(Id_{\LP(g)})\) and \(\varphi(Id_{\LP(g')})\) are two orthogonal idempotents of rank \(1\) which sum to \(I_2\). They can be simultaneously diagonalized as demanded by computing generators of their respective images.
        \end{description}

        We now prove that the algorithm is correct. First, since \(\LP(g)\) is indecomposable, by \Cref{prop:KrullSchmidtEnd} we have \(D(\LP(g) \oplus \LP(g')) \simeq M_{2}(k')\) for some finite extension \(k'/k\) if and only if \(\LP(g') \simeq \LP(g)\). Hence, our two tests do detect correctly whether \(\LP(g) \simeq \LP(g')\).
        
        Assume that \(\LP(g) \simeq \LP(g')\). Then, after conjugating by matrix \(P\) as in \Cref{line:IdemDiagonalization}, \(\varphi\) gives an isomorphism from \(S\) to the straightforward representation of \(D(\LP(g) \oplus \LP(g'))\). Then, the matrix we return corresponds to an isomorphism from \(\LP(g)\) to \(\LP(g')\).
    \end{proof}

\subsection{Algorithms for homomorphisms of lattice pairs}\label{sec:AlgosHomomorphisms}
In this section, we give algorithms related to homomorphisms of lattice pairs. All the algorithms we present rely on the computation of a pseudo-Hermite normal form of matrices with coefficients in \(A_{fi}\) and \(A_\infty\). They are therefore only polynomial time if \Cref{conj:Hermite} is assumed.

    We first give algorithms to compute kernels and images of homomorphisms of lattice pairs. Since a lattice pair is normally composed of \(A_{fi}\) and \(A_\infty\)-submodules of \(K^r\) of full rank, we are not able to give a set theoretical definition of kernels and images. Instead, we turn to a categorical approach:
    \begin{definition}\label{def:KernelImage}
        Let \(L,L'\) be lattice pairs of respective ranks \(r\) and \(r'\), and consider a homomorphism \(f\colon K^r \to K^{r'}\) from \(L\) to \(L'\) (that is, \(f(L_{fi}) \subset L'_{fi}\) and \(f(L_\infty) \subset L'_\infty\)).
        \begin{enumerate}
            \item An image of \(f\) is a pair \(I,\iota\), where \(I\) is a lattice pair of rank \(r_i = \rank(f)\) and \(\iota\colon K^{r_i} \to K^{r'}\) is an injective linear map such that \(\iota(I_{fi}) = f(L_{fi})\) and \(\iota(I_\infty) = f(L_\infty)\).
            \item A kernel of \(f\) is a pair \((\kappa,\iota)\) such that \(\kappa\) is a lattice pair of rank \(r_\kappa = r - \rank(f)\) and \(\iota\colon K^{r_\kappa} \to K^r\) is an injective linear map such that \(\iota(\kappa_{fi}) = \Ker f \cap L_{fi}\) and \(\iota(\kappa_\infty) = \Ker f \cap L_\infty\).
        \end{enumerate}
    \end{definition}

    We observe that our definitions match the definitions of kernels and images in the Abelian category of lattice pairs, and that kernels and images are therefore unique up to isomorphism. To compute such images and kernels, we adopt a similar strategy: compute the set-theoretical kernel and image and then use \Cref{thm:DimShift} below to compute a kernel and image as defined in \Cref{def:KernelImage}.

    \begin{theorem}\label{thm:DimShift}
        Let \(L\) be a lattice pair of rank \(r\), let \(S_{fi}\) be a submodule of \(L_{fi}\) of rank \(n \leq r\) and let \(S_\infty\) be a submodule of \(L_\infty\) also of rank \(n\). We further assume that \(K S_{fi} = K S_\infty\). Then, we may compute in polynomial time a lattice pair \(L'\) of rank \(n\) and a map \(f\colon L' \to L\) such that \(f(L'_{fi}) = S_{fi}\) and \(f(L'_\infty) = S_\infty\).
    \end{theorem}

    \begin{proof}
        Assume that the modules \(S_{fi}\) and \(S_\infty\) are respectively given as the images of a pseudo-matrix \(\ap,C_{fi}\) and \(C_\infty\). Then, a matrix pair \((\ap,g_{fi},g_\infty)\) and a matrix \(C \in M_{r,n}(K)\) will be a solution to the problem if
        \[C_{fi} = C g_{fi}\]
        and
        \[C_\infty = C g_\infty.\]
        We set \(g_\infty = I_n\), so the problem becomes
        \[C_{fi} = C_\infty g_{fi}.\]
        However, since the matrices \(C_{fi}\) and \(C_\infty\) both have rank \(n\) and have equal images, there exists a matrix \(g_{fi} \in GL_n(K)\) such that \(C_{fi} = C_\infty g_{fi}\) and it may be computed in polynomial time by solving a system of linear equations.
    \end{proof}

    \begin{cor}\label{cor:AlgoImage}
        Assuming \Cref{conj:Hermite}, there is an algorithm which computes the image of a homomorphism of lattice pairs in polynomial time.
    \end{cor}

    \begin{proof}
        Let \(L,L'\) be lattice pairs of respective ranks \(r\) and \(r'\) and let \(g = (\ap,g_{fi},g_\infty)\) be a matrix pair representing \(L\). Let \(f\colon L \to L'\) be a homomorphism represented by a matrix \(M \in M_{r',r}(K)\). Now, \(f(L_{fi})\) is the image of the pseudo-matrix \((\ap,Cg_{fi})\) and \(f(L_\infty)\) is \(Cg_\infty A_\infty^r\). By \Cref{prop:Hermite} (2), a pseudo-matrix of full rank spanning \(f(L_{fi})\) and a matrix of full rank spanning \(f(L_\infty)\) may be computed in polynomial time from the Hermite normal forms of pseudo-matrix \((\ap,Cg_{fi})\) over \(A_{fi}\) and matrix \(C g_\infty\) over \(A_\infty\). Then, an image of \(f\) may be computed using \Cref{thm:DimShift}.
    \end{proof}

    \begin{cor}\label{cor:AlgoKernel}
        Assuming \Cref{conj:Hermite}, there is an algorithm which computes the kernel of a homomorphism of lattice pairs in polynomial time.
    \end{cor}

    \begin{proof}
        The proof is similar to that of \Cref{cor:AlgoImage}.
    \end{proof}

    Finally, if the constant base field \(k\) is finite, we may compute a splitting of a lattice pair.

    \begin{algorithm}
        \KwIn{A matrix pair \(g\) of rank \(r\)}
        \KwOut{Matrix pairs \(g_1,\hdots,g_s\), integers \(n_1,\hdots,n_s\) and a matrix \(C \in M_r(K)\) such that the \(\LP(g_i)\) are indecomposable lattice pairs and \(C\) gives an isomorphism \(\bigoplus_{i=1}^s \LP(g_i)^n_i \simeq \LP(g)\)}
        Compute structure constants for the \(k\)-algebra \(A = \End(\LP(g))\)\; \nllabel{line:EndAlgebra}
        Compute a Wedderburn-Malcev complement \(D(A)\)\; \nllabel{line:WedderburnMalcev}
        Compute simple algebras \((S_i)_{1 \leq i \leq t}\) such that \(D(A) \simeq \bigoplus_{1 \leq i \leq t} S_i\)\;\nllabel{line:SemiSimpleSplit}
        Compute the projection maps \(p_i\colon D(A) \to S_i\)\;\nllabel{line:SemiSimpleProjs}
        \For{\(1 \leq i \leq t\)}{
            Compute \(n_i \in \N\), a finite extension \(k_i\) of \(k\) and an isomorphism \(\varphi_i\colon S_i \to M_{n_i}(k_i)\)\; \nllabel{line:MatrixIso}
            Set \(e_{ij} = (\varphi_i \circ p_i)^{-1}\left(\mathrm{Diag}(0,\hdots,1,0,\hdots,0)\right)\), with the nonzero coefficient in \(j\)-th position, for \(1 \leq j \leq n_i\)\; 
            Compute images \((g_{ij},A_{ij})\) of the endomorphisms \(e_{ij}\) of \(\LP(g)\)\; \nllabel{line:Image}
            Compute isomorphisms \(B_{ij}\colon \LP(g_{i1}) \to \LP(g_{ij})\)\; \nllabel{line:IsomsInd}
        }
        Compute the matrix \(T\) defined as the horizontal joint of the matrices \(A_{ij} B_{ij}\) as \(1 \leq i \leq t\) and \(1 \leq j \leq n_i\) are enumerated in lexicographic order\; 
        \KwRet{\(\left((g_1,\hdots,g_s),(n_1,\hdots,n_s),T\right)\)}
        \caption{Splitting a lattice pair}
        \label{algo:Splitting}
    \end{algorithm}

    \begin{theorem}\label{thm:Splitting}
        If the constant base field \(k\) is finite and \Cref{conj:Hermite} is true, \Cref{algo:Splitting} gives a correct output in polynomial time.
    \end{theorem}

    \begin{proof}
        First, we prove that every step of the algorithm makes sense and may be done in polynomial time.
            \begin{description}
                \item[\Cref{line:EndAlgebra,line:WedderburnMalcev,line:SemiSimpleSplit,line:SemiSimpleProjs,line:MatrixIso}:] This was already discussed in the proof of \Cref{lemma:AlgoIsomIndecomposable}.
                \item[\Cref{line:Image}:] can be done using \Cref{cor:AlgoImage}.
                \item[\Cref{line:IsomsInd}:] may be done using \Cref{algo:IsomIndecomposable}.
            \end{description}
            Finally, the number \(t\) of loop iterations is bounded by \(r\), the rank of \(g\). 

            Now, we prove that the output of \Cref{algo:Splitting} is correct. First, we prove that \(L \coloneq \LP(g)\) is indeed isomorphic to \(\bigoplus \LP(g_{i1})^{n_i}\). By \Cref{prop:KrullSchmidtEnd}, \[\End(L) = \bigoplus M_{n_i}(D(\End(L_i))) \oplus J(\End(L)),\] the \(S_i\) are the \(M_{n_i}(D(\End(L_i))) \simeq M_{n_i}(k_i)\) (up to reordering) and up to an automorphism of \(L\), the elements \(e_{ij}\) are the projections on a factor \(L_i\) of \(L\). An image of \(e_{ij}\) is a vector bundle \(\widetilde{L}_{ij}\) isomorphic with \(L_i\). It follows that \(L \simeq \bigoplus \LP(g_{i0})^{n_i}\).

            Then, it is easy to see that by construction, \(T\) gives an isomorphism as desired.
    \end{proof}

    \begin{cor}\label{cor:IsomFiniteField}
        If \(k\) is finite and \Cref{conj:Hermite} is true, there is an algorithm which decides if two lattice pairs are isomorphic and outputs and isomorphism if they are in polynomial time.
    \end{cor}
    \begin{proof}
        We may compute splittings for \(L\) and \(L'\) using \Cref{algo:Splitting}. Then, it is only a matter of checking that their indecomposable components are isomorphic (up to reordering) and that they appear with equal power. This may be achieved by repeated use of \Cref{algo:IsomIndecomposable}.
    \end{proof}

\section{Applications} \label{sec:Applications}
\subsection{Vector bundles on an elliptic curve}\label{sec:EllipticCurves}
In \cite{atiyah1957vector}, Atiyah gave a systematic description of the category of vector bundles on an elliptic curve over an algebraically closed field \(k\). Let \(X\) be such an elliptic curve with function field \(K\), and let \(E(r,d)\) be the set of isomorphism classes of indecomposable \(A_R\)-lattices of rank \(r\) and degree \(d\) over \(K\). In what follow, we give a succinct summary of those of his construction, rephrased in our setting of \(A_R\)-lattices, and then we give an explicit construction using lattice pairs.

\begin{definition}\label{def:AtiyahBundle}
    Let \(L\) be an \(A_R\)-lattice, let \(s = \dim_k H^0(L)\) and let \(\omega\) be a differential of \(K\). Observe that, by \Cref{prop:AbstractExt} and Serre duality, \[\Ext^1(L,R(\iota(\omega)^{-1})^s) = H^1(L^\vee \otimes R(\iota(\omega)^{-1})^s) \simeq H^0(A_R^s \otimes L)^\vee = H^0(L^s)^\vee.\]
    Upon fixing a basis of \(H^0(L)\), this extension group identifies with \(\End_k(H^0(L))\). We define the \emph{Atiyah extension} of \(L\) as the extension
    \[0 \to R(\iota(\omega^{-s})) \to L' \to L \to 0\]
    given by the identity automorphism of \(H^0(L)\).
\end{definition}

\begin{prop}
    Let \(r \in \N\). Then there exists a unique \(F_r \in E(r,0)\) such that \(\dim_k H^0(X,F_r) = 1\). For \(L \in E(r,0) \setminus \{F_r\}\), \(\dim_k H^0(L) = 0\).
\end{prop}

We let \(\Pic^0(K)\) be the group of isomorphism class of \(A_R\)-lattices of rank \(1\). We note that if \(K\) has a unique infinite places \(O\) of degree \(1\) with uniformizer \(\pi\), the elements of \(\Pic^0(K)\) are uniquely represented by \(A_R\) and the \(\LP(\mathfrak{p},1,\pi^{-1})\), where \(\mathfrak{p}\) varies over the prime ideals of \(A_{fi}\).
\begin{prop}
    Let \(r \in \N\) and \(d \in \Z\). Fix a rank one \(A_R\)-lattice \(L_1\) of degree \(1\).
    \begin{itemize}
        \item \(F_1\) is represented by \(A_R\).
        \item \(F_r\) is the Atiyah extension of \(F_{r-1}\).
        \item The map \(L \mapsto F_r \otimes L\) gives a bijection \(\Pic^0(K) \to E(r,0)\).
        \item Assume that \(d>0\), the Atiyah extension gives a bijection \(E(r,d) \to E(r+d,d)\).
        \item The map \(E \mapsto E \otimes L\) gives a bijection \(E(r,d) \to E(r,d+r)\).
        \item The map \(E \mapsto E^\vee\) gives a bijection \(E(r,d) \to E(r,-d)\).
    \end{itemize}
\end{prop}

Put together, these facts give explicit bijections \(\Pic^0(X) \to E(r,d)\) for all \(r \in \N\), \(d \in \Z\). Furthermore, the works \cite{tillmann1983unzerlegbare,arason1992indecomposable} showed that these constructions are also valid on an arbitrary perfect field \(k\) if \(E(r,d)\) now means the set of isomorphism classes of absolutely indecomposable vector bundles. Since an indecomposable vector bundle is always the trace of an absolutely indecomposable vector bundle defined over some finite extension of \(k\), this yields an algorithm for constructing any indecomposable lattice pair over an elliptic curve over a perfect field.

We also note that a generalisation to curves of genus \(1\) with no rational points was given in \cite{pumplun2004vector}.

\begin{example}\label{ex:AtiyahBundle}
    Let \(K = \F_7(x)[y]/(y^2-x^3-x)\). In this example, we construct the image of the line bundle \(\mathcal{L}(\mathfrak{p} - \infty)\) in \(E(3,2)\), where \(\mathfrak{p}\) is the prime ideal \(\langle x,y\rangle\) of \(A_{fi}\) and \(\infty\) is the divisor of the unique prime ideal of \(A_\infty\). First, we set
    \[L_0 = \LP(\mathfrak{p}^{-1},1,\pi).\]
    Following the construction of the map \(\Pic^0(X) \to E(3,2)\), we must first tensor the lattice pair \(L_0\) twice by a fixed lattice pair of degree \(1\). This will yield a lattice pair lying in \(E(1,2)\). We shall then take its Atiyah extension to get an element of \(E(3,2)\).

    We compute the tensor product \(L_1 = L_0 \otimes L_\infty^{\otimes 2}\), where \(L_\infty\) represents the line bundle of degree \(1\) \(\mathcal{L}(\infty)\). We get
    \[L_1 = \LP(\mathfrak{p}^{-1},1,\pi^{-1}).\]
    Since \(\deg(L_1) = 2\) and \(K\) has genus \(1\), it follows by the Riemann-Roch theorem that \(\dim_k(H^0(L_1)) = 2\). Applying \Cref{algo:H0}, we find that a basis for \(H^0(L_1)\) is \((1,x\pi)\).

    Now, we let \(\omega = d\pi\), and recall that this differential's divisor is the principal divisor corresponding to \(h \coloneq \frac{x^2 +3}{x^2}\), and so the corresponding lattice pair is
    \[L_\omega = \LP\left((h),1,1\right).\]
    We must now compute a répartition vector representing the element of \(H^1(L_1^\vee \otimes L_\omega^s)\) corresponding to the identity automorphism of \(H^0(L_1)\) as discussed in \Cref{def:AtiyahBundle}. This \(H^1\) group is the dual of the vector space \(H^0(L_1 \otimes L_t^s)\) by Serre duality, where \(L_t\) is the trivial lattice pair:
    \[L_t = \LP(A_{fi},1,1).\]
    Now, it is quite clear that a basis of \(H^0(L_1 \otimes L_t^s)\) is \(((1,0),(x\pi,0),(0,1),(0,x\pi))\). The space \(H^0(L_1 \otimes L_t^s)\) is identified with \(\End_k(H^0(L_1))\) by mapping a vector \((a,b)\) to the \(k\)-linear map sending \(1\) to \(a\) and \(x\pi\) to \(b\). Thus, we shall find a vector \((\alpha,\beta) \in K^2\) such that
        \begin{equation}
            \begin{cases}
                \Res_\infty(\alpha_\infty) = \Res_\infty(x\pi \beta_\infty) = 1 \\
                \Res_\infty(x\pi \alpha_\infty) = \Res_\infty(\beta_\infty) = 0.
            \end{cases}
            \label{eq:SerreDualityAtiyahEq}
        \end{equation}
        Observe that \(x\pi = \pi^{-1} + O(\pi^3)\), so we may set \(\alpha = \pi^{-1}\) and \(\beta = 1\). We have shown that the Atiyah extension of \(L_1\) is represented in \(H^1(L_1^\vee \otimes L_\omega^s)\) by the répartition vector \(\kappa = (\pi^{-1}_\infty,1_\infty)\). By \Cref{thm:ExplicitExt}, it follows that the Atiyah extension of \(L_1\) is the lattice pair
        \[L = \LP\left(((h),(h),\mathfrak{p}^{-1}),I_3, \begin{pmatrix} 1 & 0 & -\pi^{-2} \\ 0 & 1 & -\pi^{-1} \\ 0 & 0 & \pi^{-1} \end{pmatrix} \right).\]
    The determinant of \(L\) is 
    \[\det(L) = \LP\left(h^2\mathfrak{p}^{-1},1,\pi^{-1}\right).\]
    Now, the divisor of \(\mathbb{p}\) has degree \(1\) and the finite part of the divisor of \(h\) has degree \(0\), so we find that \(\deg(L) = 2\) as expected. 
     By \cite[Theorem 6]{atiyah1957vector}, we should have
    \[\det(L) \simeq L_1,\]
    and indeed one observes readily that division by \(h^2\) is such an isomorphism.

    We will now check that the lattice pair \(L\) is indeed absolutely indecomposable. By \Cref{prop:AbsolutelyIndecomposable}, we need to check that \(D(L) = \F_7\). In fact, since the rank and degree of \(L\) are coprime, we expect \(\End(L) = \F_7\) by \cite[Corollary 2.5]{oda1971vector}. We first compute \(\Endc(L,L) = L^\vee \otimes L\). Using the formulas from \Cref{def:AlgosBasics}, we find:
    \[\Endc(L) = \LP(\ap,I_3,g_\infty),\]
    with
    \[\ap = (A_{fi},A_{fi},(h\mathfrak{p})^{-1},A_{fi},A_{fi},(h\mathfrak{p}^{-1}),h\mathfrak{p},h\mathfrak{p},A_{fi})\]
    and
    \begin{align*}
        g_\infty &= \left(\begin{pmatrix} 1 & 0 & -\pi^{-2} \\ 0 & 1 & -\pi^{-1} \\ 0 & 0 & \pi^{-1} \end{pmatrix}^t\right)^{-1} \otimes \begin{pmatrix} 1 & 0 & -\pi^{-2} \\ 0 & 1 & -\pi^{-1} \\ 0 & 0 & \pi^{-1} \end{pmatrix} \\
                 &= \begin{pmatrix}1 & 0 & \frac{-x^3}{x^2+1} & 0 & 0 & 0 & 0 & 0 & 0 \\
                     0 & 1 & \frac{-xy}{x^2+1} & 0 & 0 & 0 & 0 & 0 & 0 \\
                     0 & 0 & \frac{xy}{x^2+1} & 0 & 0 & 0 & 0 & 0 & 0 \\
                     0 & 0 & 0 & 1 & 0 & \frac{-x^3}{x^2+1} & 0 & 0 & 0 \\
                     0 & 0 & 0 & 0 & 1 & \frac{-xy}{x^2+1} & 0 & 0 & 0 \\
                     0 & 0 & 0 & 0 & 0 & \frac{xy}{x^2+1} & 0 & 0 & 0 \\
                     \frac{xy}{x^2+1} & 0 & \frac{-x^4y}{(x^2+1)^2} & 1 & 0 & \frac{-x^3}{x^2+1} & \frac{y}{x^2} & 0 & \frac{-xy}{x^2+1} \\
                     0 & \frac{xy}{x^2+1} & \frac{-x^3}{x^2+1} & 0 & 1 & \frac{-xy}{x^2+1} & 0 & \frac{y}{x^2} & -1 \\
                 0 & 0 & \frac{x^3}{x^2+1} & 0 & 0 & \frac{xy}{x^2+1} & 0 & 0 & 1\end{pmatrix}
    \end{align*}

    Using the algorithm from \Cref{cor:H0}, we may compute \(\End(L) = H^0(\Endc(L))\) and find a \(1\) dimensional \(\F_7\)-vector space, whose basis element identifies with the identity matrix under our usual identification \(K^9 \simeq M_3(K)\). So, we indeed have \(\End(L) = \F_7\).

\end{example}

\subsection{Algebraic-geometry codes}\label{sec:Codes}
In \cite{savin2008algebraic}, Savin introduced a generalisation of algebraic-geometry codes to vector bundles of arbitrary rank. His construction is optimal when performed over so-called \emph{weakly-stable} vector bundle, and this motivated a line of work constructing weakly-stable vector bundles on projective curves over finite fields \cite{ballico2008vector, nakashima2010ag}. Independently, \cite{weng2018codes} gave a similar construction based on his adelic setting for vector bundles and introduced the notion of \(D\)-balanced vector bundle, where \(D\) is an effective divisor. As in \Cref{sec:EllipticCurves}, we rephrase known definitions and results in terms of \(A_R\)-lattices and give an explicit example as a lattice pair.

For what follows, we assume that \(K\) is an algebraic function field over a finite constant field \(k\).

\begin{definition}
    Let \(L\) be an \(A_R\)-lattice. The \emph{slope} of \(E\) is defined as
    \[\mu(L) \coloneq \frac{\deg(L)}{\rank(L)}.\]
\end{definition}

\begin{definition}
    An \(A_R\)-lattice \(L\) is said to be \emph{weakly-stable} if for all rank \(1\) \(A_R\)-sublattices \(L'\) of \(L\),
    \[\mu(L') \leq \mu(L).\]
\end{definition}

\begin{definition}
    Let \(D\) be an effective divisor of \(K\), and let \(L\) be an \(A_R\)-lattice. Then \(L\) is \emph{\(D\)-balanced} if \(L_P = A_P\) for all places \(P\) in the support of \(D\).
\end{definition}

\begin{prop}[\cite{savin2008algebraic}]
    Let \(r,d \in \N\). Let \(\alpha, \beta\) be the quotient and rest of the Euclidean division of \(d\) by \(r\). Let \(L'_1,L'_2\) be rank \(1\) \(A_R\)-lattices of degree \(\alpha\) and let \(L'\) be a rank \(1\) \(A_R\)-lattice of degree \(\alpha+1\). Consider the following construction:
    \begin{enumerate}
        \item \(L_1 \coloneq L'_1\).
        \item for \(2 \leq i \leq r - \beta + 1\), \(L_i\) is a non-trivial extension of \(L'_2\) by \(L_{i-1}\).
        \item for \(r - \beta + 2 \leq i \leq r\), \(L_i\) is a non-trivial extension of \(L'\) by \(L_{i-1}\).
    \end{enumerate}
    Then, \(L_r\) is a weakly-stable \(A_R\)-lattice of rank \(r\) and degree \(d\). If \(D\) is an effective divisor with support in \(M^{fi}\) such that \(L', L'_1\) and \(L'_2\) are \(D\)-balanced. Then, if the successive extensions are constructed using the algorithm from \Cref{cor:ExtensionsAlgorithm}, the lattice \(L_r\) is \(D\)-balanced.
\end{prop}

\begin{example}
    We let \(k = \F_{101}\), \(K = k(x,y)/(y^2 - x^5 - 1)\) and construct a weakly-stable vector bundle of rank \(3\) and degree \(10\) over \(K\).

    First, we set \(\pi = \frac{y}{x^3}\), a local uniformizer of \(\infty\), the unique infinite place of \(K\). We also set \(\omega = d\pi\), and we define \(\mathfrak{p}_1 = \langle x, y+1\rangle\) and \(\mathfrak{p}_2 = \langle x, y-1\rangle\), two prime ideals of \(A_{fi}\). We will build our vector bundle from the following line bundles:
        \[L = \LP(A_{fi},1,\pi^{-4}),\] 
        \[L_1 = \LP(\mathfrak{p}_1^{-3},1,1),\]
        and
        \[L_2 = \LP(\mathfrak{p_2}^{-3},1,1).\]

    We first set \(E_1 = L_1\) and compute a non-trivial extension \(E_2\) of \(L_2\) by \(E_1\). We compute a basis of \(H^0(\iota(\omega)^{-1}E_1^\vee \otimes L_2)\). This space has dimension \(1\) and is generated by \[a = \frac{-2x^2}{x^5 + 6}(y+1).\] Computing the formal series expansion of \(a\) with respect to \(\pi\), we find \[a = -2\pi + O(\pi^6).\] Therefore, the non-trivial linear form on \(H^0(\iota(\omega)^{-1}E_1^\vee \otimes L_2)\) sending \(a\) to \(1\) is represented by the infinite répartition \(b_\infty\), where \[b = \frac{-1}{2\pi^2}.\]
    We may therefore set
    \[E_2 = \LP\left((\mathfrak{p}_1^{-3},\mathfrak{p}_2^{-3}),I_2,\begin{pmatrix}1 & \frac{1}{2\pi^2} \\ 0 & 1 \end{pmatrix} \right).\]

    Our vector bundle \(E_3\) will then be constructed as a nontrivial extension of \(L\) by \(E_2\). Again, we compute a basis of \(H^0(\iota(\omega)^{-1}E_2^\vee \otimes L)\) and find:
    \[\left(\begin{pmatrix} 0 \\ \frac{x^7}{x^5 + 6} \end{pmatrix},
        \begin{pmatrix} 0 \\ \frac{x^4}{x^5 + 6}y - \frac{x^4}{x^5 + 6} \end{pmatrix},
        \begin{pmatrix} \frac{x^7}{x^5 + 6} \\ \frac{-x^8}{2(x^5 + 6} \end{pmatrix}
        \begin{pmatrix} \frac{x^4}{x^5+6}y + \frac{x^4}{x^5+6} \\ \frac{-x^5}{2(x^5+6)}y + \frac{x^5}{2(x^5 + 6)} \end{pmatrix}
\right).\]
    We let \(f_1,f_2,f_3\) and \(f_4\) be these columns.

    We shall find a vector \(\alpha \in K^2\) such that \(\theta_\omega(\cdot,\alpha_\infty)\) corresponds to the linear form \(\begin{pmatrix} 1 & 0 & 0 & 0 \end{pmatrix}\) with respect to the dual of the basis given above. That is, we must find \(\alpha_1\) and \(\alpha_2\) in \(K\) such that for all \(1 \leq i \leq 4\),
    \begin{equation}\Res_\infty\left(\sum_{j=1}^2 f_{ij} \alpha_j\right) = \begin{cases} 1 \text{ if } i = 1, \\ 0 \text{ otherwise.}\end{cases}\label{eq:ExSavin}\end{equation}
    Following \Cref{algo:H1}, we compute \(v_1 = -4\) and \(v_2 = -6\). We compute the power series expansion of the coefficients of the \(f_i\), starting at degree \(-4\) on the first row and degree \(-6\) on the second row. We get:
        \[f_1 = \begin{pmatrix} 0 \\ \pi^{-4} + O(\pi^{-2}) \end{pmatrix},\] 
        \[f_2 = \begin{pmatrix} 0 \\ \pi^{-3} + O(\pi^{-2}) \end{pmatrix},\] 
        \[f_3 = \begin{pmatrix} \pi^{-4} + O(1) \\ \frac{-1}{2} \pi^{-6} + O(\pi^{-2}) \end{pmatrix},\]
        and
        \[f_4 = \begin{pmatrix} \pi^{-3} + O(1) \\ \frac{-1}{2}\pi^{-5} + O(\pi^{-2}) \end{pmatrix}.\]

    Therefore, we must solve the linear system 
    \[X \begin{pmatrix}
        0 & 0 & 1 & 0 \\
        0 & 0 & 1 & 0 \\
        0 & 0 & 0 & 0 \\
        0 & 0 & 0 & 0 \\
        0 & 0 & \frac{-1}{2} & 0 \\
        0 & 0 & 0 & \frac{-1}{2} \\
        1 & 0 & 0 & 0 \\
        0 & 1 & 0 & 0
        \end{pmatrix} = \begin{pmatrix} 1 & 0 & 0 & 0\end{pmatrix}.\]
    An obvious solution is \(X = \begin{pmatrix} 0 & 0 & 0 & 0 & 0 & 0 & 1 & 0 \end{pmatrix}\). Bringing things together, we set \(\alpha_1 = 0\) and \(\alpha_2 = \pi^{3}\) and \Cref{eq:ExSavin} is satisfied.
    Finally, we may set
    \[E_3 = \LP\left((\mathfrak{p}_1^{-3},\mathfrak{p}_2^{-3},A_{fi}),I_3,\begin{pmatrix}1 & \frac{1}{2\pi^2} & 0 \\
            0 & 1 & \frac{-1}{\pi} \\
    0 & 0 & \frac{1}{\pi^4} \end{pmatrix}\right).\]
    
\end{example}

\bibliographystyle{acm}
\bibliography{biblio}

\begin{thebibliography}{10}

\bibitem{arason1992indecomposable}
{\sc Arason, J.~K., Elman, R., and Jacob, B.}
\newblock On indecomposable vector bundles.
\newblock {\em Comm. Algebra 20}, 5 (1992), 1323--1351.

\bibitem{atiyah1956krull}
{\sc Atiyah, M.}
\newblock On the {K}rull-{S}chmidt theorem with application to sheaves.
\newblock {\em Bull. Soc. Math. France 84\/} (1956), 307--317.

\bibitem{atiyah1957vector}
{\sc Atiyah, M.~F.}
\newblock Vector bundles over an elliptic curve.
\newblock {\em Proc. London Math. Soc. (3) 7\/} (1957), 414--452.

\bibitem{ballico2008vector}
{\sc Ballico, E.}
\newblock Vector bundles on curves over {$\Bbb F_q$} and algebraic codes.
\newblock {\em Finite Fields Appl. 14}, 4 (2008), 1101--1107.

\bibitem{biasse2017computation}
{\sc Biasse, J.-F., Fieker, C., and Hofmann, T.}
\newblock On the computation of the {HNF} of a module over the ring of integers of a number field.
\newblock {\em J. Symbolic Comput. 80\/} (2017), 581--615.

\bibitem{brenner2022deciding}
{\sc Brenner, H., and Steinbuch, J.}
\newblock Deciding stability of sheaves on curves.
\newblock {\em Internat. J. Algebra Comput. 32}, 4 (2022), 859--884.

\bibitem{cohen1993course}
{\sc Cohen, H.}
\newblock {\em A course in computational algebraic number theory}, vol.~138 of {\em Graduate Texts in Mathematics}.
\newblock Springer-Verlag, Berlin, 1993.

\bibitem{cohen1996hermite}
{\sc Cohen, H.}
\newblock Hermite and {S}mith normal form algorithms over {D}edekind domains.
\newblock {\em Math. Comp. 65}, 216 (1996), 1681--1699.

\bibitem{cohen2000advanced}
{\sc Cohen, H.}
\newblock {\em Advanced topics in computational number theory}, vol.~193 of {\em Graduate Texts in Mathematics}.
\newblock Springer-Verlag, New York, 2000.

\bibitem{degraaf1997computing}
{\sc de~Graaf, W.~A., Ivanyos, G., K\"{u}ronya, A., and R\'{o}nyai, L.}
\newblock Computing {L}evi decompositions in {L}ie algebras.
\newblock {\em Appl. Algebra Engrg. Comm. Comput. 8}, 4 (1997), 291--303.

\bibitem{decker2002sheaf}
{\sc Decker, W., and Eisenbud, D.}
\newblock Sheaf algorithms using the exterior algebra.
\newblock In {\em Computations in algebraic geometry with {M}acaulay 2}, vol.~8 of {\em Algorithms Comput. Math.} Springer, Berlin, 2002, pp.~215--249.

\bibitem{guardia2013new}
{\sc Gu\`ardia, J., Montes, J., and Nart, E.}
\newblock A new computational approach to ideal theory in number fields.
\newblock {\em Found. Comput. Math. 13}, 5 (2013), 729--762.

\bibitem{gupta2012triangular}
{\sc Gupta, S., Sarkar, S., Storjohann, A., and Valeriote, J.}
\newblock Triangular {$x$}-basis decompositions and derandomization of linear algebra algorithms over {$K[x]$}.
\newblock {\em J. Symbolic Comput. 47}, 4 (2012), 422--453.

\bibitem{harder1975cohomology}
{\sc Harder, G., and Narasimhan, M.~S.}
\newblock On the cohomology groups of moduli spaces of vector bundles on curves.
\newblock {\em Math. Ann. 212\/} (1974/75), 215--248.

\bibitem{hartshorne2013algebraic}
{\sc Hartshorne, R.}
\newblock {\em Algebraic geometry}, vol.~No. 52 of {\em Graduate Texts in Mathematics}.
\newblock Springer-Verlag, New York-Heidelberg, 1977.

\bibitem{hess2002algorithm}
{\sc Hess, F.}
\newblock An algorithm for computing {W}eierstrass points.
\newblock In {\em Algorithmic number theory ({S}ydney, 2002)}, vol.~2369 of {\em Lecture Notes in Comput. Sci.} Springer, Berlin, 2002, pp.~357--371.

\bibitem{hess2002computing}
{\sc Hess, F.}
\newblock Computing {R}iemann-{R}och spaces in algebraic function fields and related topics.
\newblock {\em J. Symbolic Comput. 33}, 4 (2002), 425--445.

\bibitem{ivanyos2018computing}
{\sc Ivanyos, G., Kutas, P., and R\'{o}nyai, L.}
\newblock Computing explicit isomorphisms with full matrix algebras over {$\Bbb F_q(x)$}.
\newblock {\em Found. Comput. Math. 18}, 2 (2018), 381--397.

\bibitem{kailath1980linear}
{\sc Kailath, T.}
\newblock {\em Linear systems}.
\newblock Prentice-Hall Information and System Sciences Series. Prentice-Hall, Inc., Englewood Cliffs, NJ, 1980.

\bibitem{lepotier1997lectures}
{\sc Le~Potier, J.}
\newblock {\em Lectures on vector bundles}, vol.~54 of {\em Cambridge Studies in Advanced Mathematics}.
\newblock Cambridge University Press, Cambridge, 1997.
\newblock Translated by A. Maciocia.

\bibitem{motsak2010graded}
{\sc Motsak, O.}
\newblock {\em Graded commutative algebra and related structures in Singular with applications}.
\newblock PhD thesis, Technischen Universität Kaiserslautern, 2010.

\bibitem{nakashima2010ag}
{\sc Nakashima, T.}
\newblock A{G} codes from vector bundles.
\newblock {\em Des. Codes Cryptogr. 57}, 1 (2010), 107--115.

\bibitem{oda1971vector}
{\sc Oda, T.}
\newblock Vector bundles on an elliptic curve.
\newblock {\em Nagoya Math. J. 43\/} (1971), 41--72.

\bibitem{pumplun2004vector}
{\sc Pumpl\"{u}n, S.}
\newblock Vector bundles and symmetric bilinear forms over curves of genus one and arbitrary index.
\newblock {\em Math. Z. 246}, 3 (2004), 563--602.

\bibitem{ronyai1990computing}
{\sc R\'{o}nyai, L.}
\newblock Computing the structure of finite algebras.
\newblock {\em J. Symbolic Comput. 9}, 3 (1990), 355--373.

\bibitem{rosen2002number}
{\sc Rosen, M.}
\newblock {\em Number theory in function fields}, vol.~210 of {\em Graduate Texts in Mathematics}.
\newblock Springer-Verlag, New York, 2002.

\bibitem{sarkar2011normalization}
{\sc Sarkar, S., and Storjohann, A.}
\newblock Normalization of row reduced matrices.
\newblock In {\em I{SSAC} 2011---{P}roceedings of the 36th {I}nternational {S}ymposium on {S}ymbolic and {A}lgebraic {C}omputation\/} (2011), ACM, New York, pp.~297--303.

\bibitem{savin2008algebraic}
{\sc Savin, V.}
\newblock Algebraic-geometric codes from vector bundles and their decoding.
\newblock arXiv preprint, 2008.
\newblock arXiv:0803.1096v1 (5 pages).

\bibitem{serre1955faisceaux}
{\sc Serre, J.-P.}
\newblock Faisceaux alg\'{e}briques coh\'{e}rents.
\newblock {\em Ann. of Math. (2) 61\/} (1955), 197--278.

\bibitem{serre1975groupes}
{\sc Serre, J.-P.}
\newblock {\em Algebraic groups and class fields}, vol.~117 of {\em Graduate Texts in Mathematics}.
\newblock Springer-Verlag, New York, 1988.
\newblock Translated from the French.

\bibitem{smith2000computing}
{\sc Smith, G.~G.}
\newblock Computing global extension modules.
\newblock vol.~29. 2000, pp.~729--746.
\newblock Symbolic computation in algebra, analysis, and geometry (Berkeley, CA, 1998).

\bibitem{stichtenoth2009algebraic}
{\sc Stichtenoth, H.}
\newblock {\em Algebraic function fields and codes}, second~ed., vol.~254 of {\em Graduate Texts in Mathematics}.
\newblock Springer-Verlag, Berlin, 2009.

\bibitem{sugahara2012adelic}
{\sc Sugahara, K.}
\newblock Adelic riemann-roch theorem on curve.
\newblock Master's thesis, Kyushu University, 2012.

\bibitem{sagemath}
{\sc {The Sage Developers}}.
\newblock {\em {S}ageMath, the {S}age {M}athematics {S}oftware {S}ystem ({V}ersion 10.3.beta8)}, 2024.
\newblock {\tt https://www.sagemath.org}.

\bibitem{tillmann1983unzerlegbare}
{\sc Tillmann, A.}
\newblock {\em Unzerlegbare Vektorb{\"u}ndel {\"u}ber algebraischen Kurven}.
\newblock PhD thesis, FernUniversit\"at, Hagen, 1983.

\bibitem{tyurin1964classification}
{\sc Tyurin, A.}
\newblock On the classification of rank 2 vector bundles over an algebraic curve of arbitrary genus.
\newblock {\em Izv. AN SSSR. Ser. Math. 28}, 1 (1964), 21--52.

\bibitem{voight2021quaternion}
{\sc Voight, J.}
\newblock {\em Quaternion algebras}, vol.~288 of {\em Graduate Texts in Mathematics}.
\newblock Springer, Cham, 2021.

\bibitem{weil1938generalisation}
{\sc Weil, A.}
\newblock Généralisation des fonctions abéliennes.
\newblock {\em J. Math. Pures. Appl. (9) 17}, 1--4 (1938), 47--87.

\bibitem{weng2018adelic}
{\sc Weng, L.}
\newblock Adelic extension classes, atiyah bundles and non-commutative codes.
\newblock arXiv preprint, 2018.
\newblock arXiv:1809.00791v1 (25 pages).

\bibitem{weng2018codes}
{\sc Weng, L.}
\newblock Codes and stability.
\newblock arXiv preprint, 2018.
\newblock arXiv:1806.04319v1 (24 pages).

\bibitem{weng2018zeta}
{\sc Weng, L.}
\newblock {\em Zeta functions of reductive groups and their zeros}.
\newblock World Scientific Publishing Co. Pte. Ltd., Hackensack, NJ, 2018.
\newblock With appendices by the author and K. Sugahara.

\end{thebibliography}
\end{document}